\documentclass{amsart}
\usepackage{amsmath,amssymb,amsthm, mathrsfs, mathtools, bm}
\usepackage{amscd,mathtools}
\usepackage[shortlabels]{enumitem}
\usepackage[utf8]{inputenc}
\usepackage[english]{babel}
\usepackage{color}
\usepackage{hyperref}
\usepackage{cite}
\usepackage{float}
\usepackage{comment}
\usepackage{mathtools}
\usepackage{csquotes}
\usepackage[normalem]{ulem}
\usepackage{caption, subcaption}
\usepackage[margin=3cm]{geometry}
\newcommand*{\rom}[1]{\expandafter\@slowromancap\romannumeral #1@}

\usepackage{tikz}
\usetikzlibrary{automata}

\usepackage{tikz-cd}
\usetikzlibrary{calc,decorations.pathreplacing}
\usetikzlibrary{arrows}
\usetikzlibrary{decorations.pathreplacing}
\usetikzlibrary{intersections}

 \newtheorem{theorem}{Theorem}[section]
  \newtheorem{proposition}[theorem]{Proposition}
  \newtheorem{corollary}[theorem]{Corollary}
  \newtheorem{lemma}[theorem]{Lemma}

    \newtheorem{question}{Question}

  \renewcommand{\theintrothm}{\Alph{introthm}}
  \renewcommand{\theintrocor}{\Alph{introcor}}

  \theoremstyle{definition}
  \newtheorem{definition}[theorem]{Definition}

  \newtheorem*{claim*}{Claim}

    \newtheorem{notation}[theorem]{Notation}

  \newtheorem*{answer*}{Answer}
  \newtheorem*{application*}{Application}

  \newtheorem{remark}[theorem]{Remark}
  \newtheorem*{remark*}{Remark}

\newcommand{\Cone}{\text{Cone}}
\newcommand{\Pone}{\text{Pcone}}
\DeclarePairedDelimiterX{\Norm}[1]{\lVert}{\rVert}{#1}
\theoremstyle{definition}

 \newcommand{\sA}{{\sf A}}   
  \newcommand{\sB}{{\sf B}}   
  \newcommand{\sC}{{\sf C}}   
  \newcommand{\sD}{{\sf D}}     
  \newcommand{\sE}{{\sf E}}   
  \newcommand{\sF}{{\sf F}}   
  \newcommand{\sG}{{\sf G}}   
  \newcommand{\sH}{{\sf H}}     
  \newcommand{\sI}{{\sf I}}   
  \newcommand{\sJ}{{\sf J}}   
  \newcommand{\sK}{{\sf K}}   
  \newcommand{\sL}{{\sf L}}     
  \newcommand{\sM}{{\sf M}}   
  \newcommand{\sN}{{\sf N}}   
  \newcommand{\sO}{{\sf O}}   
  \newcommand{\sP}{{\sf P}}       
  \newcommand{\sQ}{{\sf Q}}   
  \newcommand{\sR}{{\sf R}}   
  \newcommand{\sS}{{\sf S}}     
  \newcommand{\sT}{{\sf T}}   
  \newcommand{\sU}{{\sf U}}   
  \newcommand{\sV}{{\sf V}}   
  \newcommand{\sW}{{\sf W}}     
  \newcommand{\sX}{{\sf X}}   
  \newcommand{\sY}{{\sf Y}}   
  \newcommand{\sZ}{{\sf Z}}

  \renewcommand{\aa}{{\sf a}}   
  \newcommand{\bb}{{\sf b}}   
  \newcommand{\cc}{{\sf c}}   
  \newcommand{\dd}{{\sf d}}     
    \newcommand{\E}{{\mathbb{e}}}   
  \newcommand{\fff}{{\sf f}}   
  \newcommand{\hh}{{\sf h}}     
  \newcommand{\ii}{{\sf Ii}}   
  \newcommand{\jj}{{\sf j}}   
  \newcommand{\kk}{{\sf k}}   
  \renewcommand{\ll}{{\sf l}}     
  \newcommand{\mm}{{\sf m}}   
  \newcommand{\nn}{{\sf n}}   
  \newcommand{\oo}{{\sf o}}   
  \newcommand{\pp}{{\sf p}}       
  \newcommand{\qq}{{\sf q}}   
  \newcommand{\rr}{{\sf r}}  
    \newcommand{\RR}{{\mathbb{R}}}    
  \renewcommand{\ss}{{\sf s}}     
  \renewcommand{\tt}{{\sf t}}   
  \newcommand{\uu}{{\sf u}}   
  \newcommand{\vv}{{\sf v}}   
  \newcommand{\ww}{{\sf w}}     
  \newcommand{\xx}{{\sf x}}   
  \newcommand{\yy}{{\sf y}}   
  \newcommand{\zz}{{\sf z}}   
  
    \newcommand{\bfa}{{\textbf{a}}} 
  \newcommand{\bfb}{{\textbf{b}}} 
    \newcommand{\bfc}{{\textbf{c}}} 
   
  \newcommand{\gothic}{\mathfrak}
  \newcommand{\go}{{\gothic o}}
  \newcommand{\calA}{\mathcal{A}}
  \newcommand{\calB}{\mathcal{B}}
  \newcommand{\calC}{\mathcal{C}}
  \newcommand{\calD}{\mathcal{D}}
  \newcommand{\calE}{\mathcal{E}}
  \newcommand{\calF}{\mathcal{F}}
  \newcommand{\calG}{\mathcal{G}}
  \newcommand{\calH}{\mathcal{H}}
  \newcommand{\calI}{\mathcal{I}}
  \newcommand{\calJ}{\mathcal{J}}
  \newcommand{\calK}{\mathcal{K}}
  \newcommand{\calL}{\mathcal{L}}
  \newcommand{\calM}{\mathcal{M}}
  \newcommand{\calN}{\mathcal{N}}
  \newcommand{\calO}{\mathcal{O}}
  \newcommand{\calP}{\mathcal{P}}
  \newcommand{\calQ}{\mathcal{Q}}
  \newcommand{\calR}{\mathcal{R}}
  \newcommand{\calS}{\mathcal{S}}
  \newcommand{\calT}{\mathcal{T}}
  \newcommand{\calU}{\mathcal{U}}
  \newcommand{\calV}{\mathcal{V}}
  \newcommand{\calW}{\mathcal{W}}
  \newcommand{\calX}{\mathcal{X}}
  \newcommand{\calY}{\mathcal{Y}}
  \newcommand{\calZ}{\mathcal{Z}}
    \newcommand{\blue}{\color{blue}}
    \DeclareMathOperator{\op}{op}
\DeclareMathOperator{\Ima}{Im}
\DeclareMathOperator{\diam}{diam}

\renewcommand{\o}{\circ}
\newcommand{\wt}{\widetilde}
\newcommand{\R}{\mathbb{R}}
\newcommand{\Z}{\mathbb{Z}}
\newcommand{\N}{\mathbb{N}}
\renewcommand{\H}{\mathbb{H}}
\newcommand{\Hom}{\mathrm{Hom}}
\newcommand{\Ker}{\mathrm{Ker}}
\newcommand{\im}{\mathrm{Im}}
\newcommand{\s}{\sigma}
\newcommand{\ra}{\rightarrow}
\newcommand{\Ra}{\Rightarrow}
\newcommand{\LRa}{\Leftrightarrow}
\newcommand{\cu}{\subseteq}
\newcommand{\g}{\gamma}
\newcommand{\G}{\Gamma}
\newcommand{\td}{\tilde}
\newcommand{\mbb}{\mathbb}
\newcommand{\mc}{\mathcal}
\newcommand{\mf}{\mathfrak}
\newcommand{\x}{\times}
\newcommand{\eps}{\epsilon}
\newcommand{\Om}{\Omega}
\newcommand{\om}{\omega}
\newcommand{\Aut}{\mathrm{Aut}}
\newcommand{\acts}{\curvearrowright}
\newcommand{\wh}{\widehat}
\newcommand{\mscr}{\mathscr}
\newcommand{\lcrt}{\mathrm{lcrt}}
\newcommand{\lcr}{\mathrm{lcr}}
\newcommand{\crt}{\mathrm{crt}}
\newcommand{\Cr}{\mathrm{cr}}
\newcommand{\Sym}{\mathrm{Sym}}
\newcommand{\lk}{\mathrm{lk}}
\newcommand{\hull}{\mathrm{Hull}}
\renewcommand{\ll}{\llbracket}
\newcommand{\CAT}{{\rm CAT(0)}}
\newcommand{\B}{\mbb{B}}
\newcommand{\X}{\mbb{X}}
\newcommand{\Hyp}{\mathcal{HYP}}
\newcommand{\FG}{\mathcal{FG}}
\newcommand{\FQ}{\mathcal{FQ}}
\newcommand{\Vis}{\mathcal{VIS}}
\newcommand{\Con}{\mathcal{CONE}}

\newcommand{\SM}{\mathcal{SM}}
\newcommand{\supp}{\mathrm{supp}}

\newcommand{\bma}{\bm{b}} 
\newcommand{\bmb}{\bm{\beta}} 

\newcommand{\orange}[1]{\textcolor{orange}{#1}}
\newcommand{\cyan}[1]{\textcolor{cyan}{#1}}
\newcommand{\red}[1]{\textcolor{red}{#1}}

\newcommand{\xspace}{\,}

\newcommand{\specialcell}[2][c]{%
  \begin{tabular}[#1]{@{}c@{}}#2\end{tabular}}

  \newcommand{\ST}{\mathbin{\Big|}} 
\newcommand{\pka}{\partial_{\kappa}} 
\newcommand{\sg}{\mathfrak{g}} 

\newcommand{\ou}{%
  \mathrel{%
    \vcenter{\offinterlineskip
      \ialign{##\cr$+$\cr\noalign{\kern-1.5pt}$-$\cr}%
    }%
  }%
}

\newcommand{\myGlobalTransformation}[2]
{
    \pgftransformcm{0.2}{0.8}{0}{1}{\pgfpoint{#1}{#2}}
}

\begin{document}

\title{Sublinearly Morse boundaries from the viewpoint of combinatorics}


 \author{Merlin Incerti-Medici}
\address{Institut f\"ur Mathematik, ETH Z\"urich, Switzerland}
\email{merlin.incertimedici@math.ethz.ch}
 
  \author   {Abdul Zalloum}
 \address{Department of Mathematics, Queen's University, Kingston, ON }
 \email{az32@queensu.ca}

 

\begin{abstract}

We prove that the sublinearly Morse boundary of every known cubulated group continuously injects in the Gromov boundary of a certain hyperbolic graph. We also show that for all CAT(0) cube complexes, convergence to sublinearly Morse geodesic rays has a simple combinatorial description using the hyperplanes crossed by such sequences. As an application of this combinatorial description, we show that a certain subspace of the Roller boundary continously surjects on the subspace of the visual boundary consisting of sublinearly Morse geodesic rays.

\end{abstract}

\maketitle

\section{Introduction}


\subsection{Main Results}

Boundaries at infinity are a common tool in the study of large scale geometric properties. When a group acts geometrically on a metric space and the space satisfies some suitable curvature conditions, the boundary carries rich topological, dynamic, metric, quasi-conformal, measure-theoretic and algebraic structures that allow us to study the group we started with. 
A particularly fruitful instance of this is the case of a Gromov-hyperbolic metric space and its visual boundary. Since any quasi-isometry between hyperbolic metric spaces induces a homeomorphism on the visual boundary, we can define the visual boundary of a hyperbolic group as a topological space up to homeomorphism (see \cite{Gromov1987}). This is no longer true, when the group is acting on a non-positively curved space, e.\,g.\,a $\CAT$ space. In \cite{CrokeKleiner}, Croke and Kleiner provided an example of a group acting geometrically on two different quasi-isometric $\CAT$ spaces which have non-homeomorphic visual boundaries. So we cannot associate the visual boundary as a topological space up to homeomorphism to a group.

Several attempts have recently been made to rectify this issue, most recent of which is work of Qing, Rafi and Tiozzo where they introduce the sublinearly Morse boundary and show that it's a metrizable quasi-isometry invariant of any CAT(0) space\cite{QRT19}. Given a sublinear function $\kappa$, a geodesic ray $b$ is said to be $\kappa$-Morse if there exists a function $\mm_b:\mathbb{R}^{+} \times \mathbb{R}^{+} \rightarrow \mathbb{R}^+$ such that for any $(\qq, \sQ)$-quasi geodesic $\alpha$ with end points on $b$, if $p \in \alpha,$ we have $d(p, b) \leq \mm_b(\qq,\sQ) \kappa(||p||)$, where $||p||=d(p, b(0)).$ Elements of the $\kappa$-Morse boundary are in a one to one correspondence with $\kappa$-Morse geodesic rays starting at a fixed point $\go \in X.$ The advantage of the $\kappa$-Morse boundary over other hyperbolic-like boundaries is that it's generic with respect to various measures (see Subsection \ref{subsec: history} on the history of such boundaries). Therefore, statements about the  $\kappa$-Morse boundary could be regarded as statements about almost every (quasi)-geodesic ray in $X$. Since the $\kappa$-Morse boundary of a CAT(0) space $X$ is a quasi-isometry invariant, we may denote it by $\partial_\kappa G$ whenever $G$ is a group acting geometrically on $X$. Our first result is that if a group $G$ admits a geometric action on a CAT(0) cube complex with a factor system, then the $\kappa$-Morse boundary of $G$ continuously injects in the Gromov boundary of a certain geodesic hyperbolic graph $\Gamma.$ We remark that all known proper cocompact cube complexes admit a factor system, see \cite{HagenSusse} by Hagen and Susse.

\begin{theorem}\label{introthm:cont injection}
 
  Let $G$ be a group acting geometrically on a \CAT \, cube complex $X$ with a factor system. There exist a $\delta$-hyperbolic graph $\Gamma$, depending only on $X$, and a projection map $p:X \rightarrow \Gamma$ such that for any sublinear function $\kappa$ the following holds:
  
  \begin{enumerate}
      \item  Every $\kappa$-Morse geodesic ray in $X$ projects to an infinite unparameterized quasi-geodesic in $\Gamma.$
      \item  The $\kappa$-Morse boundary of $G$ denoted $\partial_\kappa G$ continuously injects in the Gromov boundary $\partial \Gamma.$ 
    \end{enumerate}
    When $G$ is a right-angled Artin group, the graph $\Gamma$ is the contact graph of its corresponding Salvetti complex.
 \end{theorem}

We remark that in the special case of the Morse boundary ($\kappa=1$), the above theorem is known to hold in the setting of hierarchically hyperbolic spaces due to Abbott, Behrstock, and  Durham \cite{Abbot}. More precisely, they show that if $\mathcal{X}$ is a hierarchically hyperbolic space with mild assumptions, then there exists a $\delta$-hyperbolic space $Y$ such that every Morse geodesic ray in $\mathcal{X}$ projects to an infinite unparameterized quasi-geodesic in $Y$ defining a continuous injection between the Morse boundary of $\mathcal{X}$ and the Gromov boundary of $Y.$ Therefore, it's natural to wonder if the same conclusion of Theorem \ref{introthm:cont injection} holds in the settings of hierarchically hyperbolic spaces. Using work of Rafi and Verberne \cite{RafiVerberne}, we show that the answer to this questions is negative.

\begin{proposition}\label{introprop: Failure in HHG} There exists a hierarchically hyperbolic space $\mathcal{X}$ such that the following holds. If $Y$ is the maximal $\delta$-hyperbolic space associated to $\mathcal{X}$, and $p:\mathcal{X} \rightarrow Y$ is the projection map to $Y$, then all Morse geodesic rays in $\mathcal{X}$ project to infinite unparameterized quasi-geodesics in $Y$. Nonetheless, there exists a $\kappa$-Morse geodesic ray in $\mathcal{X}$ which projects to an infinite diameter set that is not an unparameterized quasi-geodesic in $Y$.

\end{proposition}


Denote the boundary of CAT(0) space $X$ as a set by $\partial X$. On the other hand, when we equip this set $\partial X$ with the standard cone topology, we denote it by $\partial_{\infty} X$ and refer to it as the visual boundary. When considering a $\CAT$ cube complex $X$, one can define the following topology on the set $\partial X$: Fix a vertex $\go \in X$ as a base point and let $h_1, \dots, h_n$ be distinct hyperplanes in $X$. 

\begin{equation*}
\begin{split}
V_{\go, h_1, \dots, h_n} := \{ \xi \in \partial X | \text{ The unique} & \text{ geodesic representative of $\xi$ based}\\
& \text{ at $\go$ crosses the hyperplanes $h_1, \dots, h_n$} \}.
\end{split}
\end{equation*}

The collection $B=\{ V_{\go, h_1, \dots, h_n}|  n \in \mathbb{N}, \,h_1,h_2,..,h_n \text{ are hyperplanes} \}$ forms the basis of a topology which we denote $\Hyp$.

It is worth noting that the $\Hyp$-topology is usually very different from the cone topology. For example, in $\mathbb{R}^2,$ the cone topology on $\partial \mathbb{R}^2$ gives a circle, however, the $\Hyp$-topology on $\partial \mathbb{R}^2$ is not even Hausdorff.

While the  cone topology is different from the $\Hyp$-topology on the set $\partial X$, our second main result states that the two topologies agree when restricted to the subset of $\partial X$ consisting of all $\kappa$-Morse geodesic rays, we denote this subset by $\partial^\kappa X \subseteq \partial X.$

\begin{theorem}\label{introthm: HYP=sublinear}
Let $X$ be a finite dimensional $\CAT$ cube complex and let $\kappa$ be a sublinear function. The restrictions of the cone topology and $\Hyp$ to the subset $\partial^{\kappa} X \subset \partial X$ are equal.
\end{theorem}

The combinatorial nature of $\Hyp$ enables us to understand some connections between the $\kappa$-Morse boundary and the Roller boundary $\partial_R X$ (see Section \ref{subsec:RollerBoundary} for a precise definition). The Roller boundary can be thought of as a combinatorial counterpart to the visual boundary. Fixing a vertex $\go$ as base point in $X$, every point $x$ in the Roller boundary can be represented by a combinatorial geodesic ray $\alpha_x$ in the $1$-skeleton of $X$ that starts at $\go$. In general, such a combinatorial geodesic ray doesn't uniquely determine a point in the visual boundary of $X.$ However, we show that when the combinatorial geodesic ray $\alpha_x$ is $\kappa$-Morse, every other associated combinatorial geodesic ray $\beta_x$ is also $\kappa$-Morse, and that $x$ determines a unique point in the visual boundary. We denote the set of points in the Roller boundary that can be represented by a $\kappa$-Morse combinatorial geodesic ray by $\partial_R^{\kappa} X$. We also denote the subset of the visual boundary consisting of points that can be represented by $\kappa$-Morse geodesic rays by $\partial^{\kappa}_\infty X$ when equipped with the subspace topology.

\begin{theorem}\label{introthm: Roller map}
Let $X$ be a finite dimensional $\CAT$ cube complex. There exists a well-defined continuous surjective map $\Phi : \partial_R^{\kappa} X \rightarrow \partial^\kappa_{\infty} X$. For all $x, y \in \partial_R^{\kappa} X$, we have that $\Phi(x) = \Phi(y)$ if and only if $[x \vert y]_{\go} = \infty$. (See Definition \ref{def:CombinatorialGromovProduct} for a definition of $[\cdot \vert \cdot]_{\go}$).

Furthermore, the induced quotient map $\overline{\Phi} : \overline{\partial_R^{\kappa} X} \rightarrow \partial_\infty^{\kappa} X$ is a homeomorphism.
\end{theorem}

We remark that the above theorem provides an extension of Theorem D of Beyrer, Fioravanti in \cite{BeyFio} and Theorem 1.1 of Zalloum in \cite{Zal18}. More precisely, Theorem D in \cite{BeyFio} shows Theorem \ref{introthm:cont injection} for the special case where $\kappa=1$ (the Morse boundary of a CAT(0) cube complex). On the other hand, Theorem 1.1 in \cite{Zal18} shows that for any proper geodesic metric space $X$, an appropriate subspace of the horofunctions boundary of $X$ continuously surjects on the Morse boundary of $X$. When $X$ is a CAT(0) cube complex, the horofunctions boundary with respect to the $l^1$-metric of $X$ is exactly the Roller boundary (Proposition 6.20 of \cite{Ferns2018} by Fernós, Lécureux and Mathéus).

\subsection{Sketch of the proof of Theorem \ref{introthm:cont injection}}

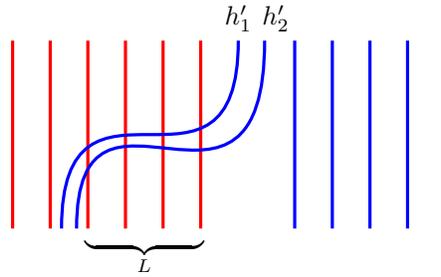
\begin{figure}

$$
\begin{tikzpicture}[scale=.50]

\draw[ very thick,red] (0,0) -- ++(0,5);

\draw[ very thick,red] (-1,0) -- ++(0,5);
\draw[ very thick,red] (-2,0) -- ++(0,5);
\draw[ very thick,red] (-3,0) -- ++(0,5);
\draw[ very thick,red] (-4,0) -- ++(0,5);

\draw[ very thick,red] (1,0) -- ++(0,5);
\draw[ very thick,blue] (3.5,0) -- ++(0,5);
\draw[ very thick,blue] (4.5,0) -- ++(0,5);
\draw[ very thick,blue] (5.5,0) -- ++(0,5);
\draw[ very thick,blue] (6.5,0) -- ++(0,5);
\draw[ very thick,blue] (7.5,0) -- ++(0,5);

\draw[very thick,blue] (2,5) .. controls ++(0,-5) and ++(0,5) ..(-2.7,0);
\draw[very thick,blue] (2.7,5) .. controls ++(0,-6) and ++(0,5) ..(-2.3,0);

\node[below] at (-.5,0) {$\underbrace{\hspace{45pt}}_L$};

\node[above] at (2,5) {$h'_{1}$};
\node[above] at (3,5) {$h'_{2}$};

\end{tikzpicture}
$$

\caption{The first two blue hyperplanes $h_1', h_2'$ can cross at most $L$ of the red hyperplanes.}\label{fig: bounded crossings}
\label{fig:ProjectiontoUnparametrisedGeodesics}

\end{figure}

We fix a finite dimensional $\CAT$ cube complex $X$.  A collection of three disjoint hyperplanes in $X$ is said to form a \emph{facing triple} if none of them separates the other two. Two disjoint hyperplanes $h_1,h_2$ in $X$ are said to be \emph{$k$-well-separated} if the number of hyperplanes meeting them both and containing no facing triple is bounded above by $k$ (for a detailed comparison between the notion of $k$-well-separation versus the Charney and Sultan's $k$-separation notion, see Section 1.2 of \cite{Murray-Qing-Zalloum}). In \cite{HypInCube}, starting with a $\CAT$ cube complex $X$, Anthony Genevois constructs a (non-geodesic) hyperbolic metric space $(Y_k,d_k)$ where $Y_k$ is the collection of vertices $X^{(0)}$ and $d_k(x,y)$ is the cardinality of the maximal collection of $k$-well-separated hyperplanes separating $x,y$, in fact, he shows that when $k=0,$ the space $Y_k$ is quasi-isometric to the contact graph of Hagen \cite{Hagen2014}. An essential result due to Genevois in \cite{HypInCube} is that for a cocompact CAT(0) cube complex $X$ with a factor system, there exists an $L$ depending only on $X$, such that if two hyperplanes are $k$-well-separated for some $k,$ then $k \leq L.$ That is to say, two hyperplanes are either not well-separated or $L$-well-separated. We will refer to this constant $L$ as the \emph{separation constant} of $X$. Let $G$ be a group admitting a geometric action on a CAT(0) cube complex $X$ with a factor system, and let $L$ be the separation constant, our steps for the proof are as follows:
 
 \begin{enumerate}

     \item Denote $\partial_\kappa X$ the quasi-isometry invariant $\kappa$-Morse boundary of $X$ given in \cite{QRT19}, and $\partial_\infty^\kappa X$ the subspace of the visual boundary of $X$ consisting of $\kappa$-Morse geodesic rays with the subspace topology. We first show that the map
     
     $$\iota_1:\partial_\kappa X \rightarrow \partial_\infty^\kappa X,$$ is continuous. That is to say, the quasi-isometry invariant topology of Qing, Rafi and Tiozzo \cite{QRT19} is finer than the subspace topology induced on the collection of $\kappa$-Morse geodesic rays. This is Lemma \ref{lem:injection}.
      \item  We construct a graph $(\Gamma, d_{\Gamma})$ which is bilipschitz equivalent to $(Y_L,d_L).$ The graph is built as follows: the vertices are the same as the vertices of $Y_L$ which are simply the vertices of $X^{(0)}$, and we connect two vertices $x,y$ with an edge if the number of $L$-well-separated hyperplanes separating them is at most $10L+4.$ This is Lemma \ref{lem:updating to a graph}.
     
     \item We show that every (CAT(0) or combinatorial) geodesic of $X$ projects to an unparameterized quasi-geodesic in $(Y_L,d_L)$ and hence to $(\Gamma,d_{\Gamma}).$ Since the argument for this is short and combinatorial, we sketch it here. Let $c$ be a geodesic in $X$ and let $x,y,z$ be the points $c(t_1),c(t_2),c(t_3)$ for $t_1 < t_2 <t_3.$ Suppose that $d_k(x,y)=L_1$ and that $d_k(y,z)=L_2$. That is to say, if $\mathcal{H}_1, \mathcal{H}_2$ are the maximal collections of $L$-well-separated hyperplanes separating $x,y$ and $y,z$ respectively, then $|\mathcal{H}_1|=L_1$ and $|\mathcal{H}_2|=L_2$. We will think of the sets $\mathcal{H}_1, \mathcal{H}_2$ as ordered sets based on the order they are crossed by $c.$ We will show that $d_L(x,z) \geq d_L(x,y)+d_L(y,z)-(L+2).$ Let $h_{\text{last}}$ denote the last hyperplane in $\mathcal{H}_1$ and $h'_{1}, h'_{2}$ denote the first and second hyperplanes of $\mathcal{H}_2.$ If $h_{\text{last}}$ and $h'_{2}$ are disjoint, there is nothing to prove. Otherwise, if $h'_{2}$ crosses some collection of $\mathcal{H}_1$, then $h'_{1}$ must also cross the same collection as  $h'_{1}$,  $h'_{2}$ are disjoint. However, since the collection $\mathcal{H}_2$ consists of $L$-well-separated hyperplanes, there can be at most $L$ hyperplanes crossed by both  $h'_{1}, h'_{2}$ (see Figure \ref{fig: bounded crossings}). Denote such hyperplanes $h_1,h_2,\dots,h_L$. The set $\mathcal{H}_1 \cup \mathcal{H}_2 \setminus\{h_1,h_2, \dots, h_L,h'_{1}, h'_{2}\}$ is a collection of $L$-well-separated hyperplanes separating $x,z$ whose cardinality is at least $L_1+L_2-(L+2)$. This shows that a geodesic $c$ in $X$ projects to an unparameterized quasi-geodesic in $Y_L,$ and by step 2, the geodesic $c$ projects to an unparameterized quasi-geodesic in $\Gamma$. This is Lemma \ref{lem:ProjectingtoQuasiGeodesics}.
     
     \item In \cite{Murray-Qing-Zalloum}, Murray, Qing and Zalloum show that a geodesic ray is $\kappa$-Morse if and only if there exists a constant $\nn$ such that $c$ crosses an infinite sequence of hyperplanes $h_i$ at points $c(t_i)$ with $d(t_i,t_{i+1}) \leq \nn \kappa(t_{i+1})$ and $h_i,h_{i+1}$ are $L$-well-separated. Hence, by part (3), a $\kappa$-Morse geodesic ray projects to an infinite unparameterized quasi-geodesic in $\Gamma.$ This defines a map $\iota_2: \partial ^\kappa_\infty X \rightarrow \partial \Gamma,$ we show this map is continuous and injective. Continuity and injectivity of the map are Lemma \ref{lem: continuous map gnrl settings} and Lemma \ref{lem: injectivity of the map} respectively.
     
     \item Composing the maps from step (1) and step (4), we get a continuous injection $\iota=\iota_2 \iota_1$
     
     $$\iota:\partial_\kappa X \rightarrow \partial \Gamma.$$ Since $\partial_\kappa X$ is independent of the space $X$ on which $G$ admits a geometric action, it may be denoted by $\partial_\kappa G.$
 \end{enumerate} 
 
\subsection{History of hyperbolic-like boundaries and cube complexes} \label{subsec: history}
The $\kappa$-Morse boundary of Qing, Rafi and Tiozzo was preceded by a few other hyperbolic-like boundaries which are also invariant under quasi-isometries. Charney and Sultan introduced the Morse boundary for $\CAT$ spaces whose homeomorphism class is invariant under quasi-isometry (see \cite{CharneySultan}). The Morse boundary was later generalized by Cordes to all proper geodesic metric spaces (see \cite{Cordes}). In either case, this boundary consists of equivalence classes of geodesic rays that have similar properties to geodesic rays in hyperbolic spaces, indicating that these are the `hyperbolic' directions in the space under consideration.

The Morse boundary of Charney-Sultan and Cordes is equipped with a quasi-isometry invariant topology by using a direct limit construction. This direct limit construction was later refined to a direct limit of Gromov hyperbolic boundaries by Cordes and Hume (see \cite{CordesHume}). Group actions on their Morse boundaries have similar dynamical properties to that of the action of a hyperbolic group on it's Gromov's boundary (see \cite{Liu19} by Liu). However, this topology is not first countable in general. Concretely, this topology is not first countable for the group $\mathbb{Z}^2 * \mathbb{Z}$, an example also known as the `tree of flats' (see \cite{Murray2019} by Murray). Cashen and Mackay introduced a refined topology for the Morse boundary (see \cite{CashenMackay}) and were able to show that the Morse boundary, with respect to their topology is a metrizable quasi-isometry invariant. The $\kappa$-Morse boundary was introduced by
Qing, Rafi and Tiozzo in \cite{QRT19} and sought to rectify some of the shortcomings found in the Morse boundary, for example, unlike the Morse boundary of Charney-Sultan, Cordes and Cashen-Mackay, the $\kappa$-Morse boundary serves as a topological model of the Poisson boundary for right-angled Artin groups (see \cite{Qing2019} and \cite{QRT19}). Furthermore, it was recently announced by Gekhtman, Qing and Rafi \cite{GQR} that $\kappa$-Morse geodesics are generic in rank-1 CAT(0) spaces in various natural measures.

In this article, we restrict our attention to the $\kappa$-Morse boundary of $\CAT$ cube complexes. Cube complexes were introduced by Gromov in \cite{Gromov1987} and have become a central object in geometric group theory over the last decade due to their fruitfulness in solving problems in group theory and low-dimensional topology and due to the fact that many interesting groups are cubulable, i.\,e.\,they act properly and cocompactly on a $\CAT$ cube complex. The class of groups that are cubulable includes Right-angled Artin groups, hyperbolic $3$-manifold groups (\cite{BergersonWise}), most non-geometric $3$-manifold groups (\cite{PiotrWise14}, \cite{HagenPiotr}, \cite{PiotrWise18}), small cancelation groups (\cite{Wise04}) and many others. Cubulated groups played a key role in Agol's and Wise's proof of the virtual Haken and virtual fibered conjecture (\cite{Wise11}, \cite{Agol13}). The additional structure of $\CAT$ cube complexes allows us to study sublinearly Morse boundaries via more combinatorial methods.

\subsection*{Outline of the paper} Section \ref{sec:Preliminaries} contains some background material on CAT(0) spaces, cube complexes, sublinearly Morse boundaries, Roller boundaries and the various topologies these boundaries are equipped with. In Section \ref{sec:CombinatorialGromovProducts} we show that any two combinatorial geodesic rays which are crossed by an infinite sequence of hyperplanes must fellow travel, provided that at least one of them is $\kappa$-Morse. This will be essential to establishing well-definness of the map in Theorem \ref{introthm: Roller map}. Section \ref{sec:ComparingTopologies} proves Theorem \ref{introthm: HYP=sublinear} giving a combinatorial description of convergence to sublinearly Morse geodesic rays, this will be used in proving continuity of the map in Theorem \ref{introthm: Roller map}. In Section \ref{sec:RollerContinuity}, we prove Theorem \ref{introthm: Roller map}. In Section \ref{sec:EmbeddingintoGromovboundary}, we describe the hyperbolic graph $\Gamma$ and we prove Theorem \ref{introthm:cont injection}. In Section \ref{sec: remark on hhg} we prove Proposition \ref{introprop: Failure in HHG}. We remark that Sections \ref{sec:EmbeddingintoGromovboundary} and \ref{sec: remark on hhg} are independent of other sections and can be read separately.

\subsection*{Acknowledgement} The authors would like to thank Carolyn Abbott, Ruth Charney, Matthew Durham, Elia Fioravanti, Anthony Genevois, Mark Hagen, Qing Liu, Kasra Rafi and Jacob Russell for fruitful discussions.

\section{Preliminaries} \label{sec:Preliminaries}




\subsection{CAT(0) geometry} \label{subsec:CAT(0)Geometry}


\begin{definition}[Quasi-geodesics]\label{Def:Quasi-Geodesic} A \emph{geodesic ray} in a metric space $X$ is an isometric embedding $b : [0, \infty) \to X$. We fix a base-point $\go \in X$ and always assume 
that $b(0) = \go$, that is, a geodesic ray is always assumed to start from 
this fixed base-point. A \emph{quasi-geodesic ray} is a continuous quasi-isometric 
embedding $\beta : [0, \infty) \to X$ starting from $\go$. 
The additional assumption that quasi-geodesics are continuous is not necessary,
but is added to make the proofs simpler. 

\end{definition}

Let $(X,d_X)$ be a metric space. A metric space is called proper if closed balls are compact. It is called geodesic if any two points $x, y \in X$ can be connected by a geodesic segment. A proper, geodesic metric space $(X, d_{X})$ is $\CAT$ if geodesic triangles in $X$ are at 
least as thin as triangles in Euclidean space with the same side lengths. To be precise, for any 
given geodesic triangle $\triangle pqr$, consider the up to isometry unique triangle 
$\triangle \overline p \overline q \overline r$ in the Euclidean plane with the same side 
lengths. For any pair of points $x, y$ on the triangle, for instance on edges $[p,q]$ and $[p, r]$ of the 
triangle $\triangle pqr$, if we choose points $\overline x$ and $\overline y$  on 
edges $[\overline p, \overline q]$ and $[\overline p, \overline r]$ of 
the triangle $\triangle \overline p \overline q \overline r$ so that 
$d_X(p,x) = d_{\mathbb{E}^2}(\overline p, \overline x)$ and 
$d_X(p,y) = d_{\mathbb{E}^2}(\overline p, \overline y)$, then
\[ 
d_{X} (x, y) \leq d_{\mathbb{E}^{2}}(\overline x, \overline y).
\] 

For the remainder of the paper, we assume $(X, d)$ is a proper $\CAT$ space. Here, we list some properties of proper $\CAT$ spaces that are needed later (see 
 \cite{CAT(0)reference}). 

\begin{lemma} \label{Lem:CAT} 
 A proper \CAT space $X$ has the following properties:
\begin{enumerate}
\item It is uniquely geodesic, that is, for any two points $x, y$ in $X$, 
there exists exactly one geodesic connecting them. 

\item The nearest-point projection from a point $x$ to a geodesic line $b$ 
is a unique point denoted $x_b$. The closest-point projection map
\[
\pi_b : X \to b
\]
is Lipschitz. 
\item Convexity: For any convex set $Z \in X$, the distance function $f:X \rightarrow \mathbb{R}^+$ given by $f(x)=d(x,Z)$ is convex.
\end{enumerate}
\end{lemma}

\subsection{The cone topology on the visual boundary}

 \begin{definition}[space of geodesic rays] \label{def: visual topology}
As a set, the \emph{boundary} of $X$, denoted by $\partial X$ is defined to be the collection of equivalence classes of all infinite geodesic rays. Let $b$ and $c$ be two infinite geodesic rays, not necessarily starting at the same point. We define an equivalence relation as follows:  $b$ and $c$ are in the same equivalence class,  if and only if there exists some $\nn \geq 0$ such that $d(b(t), c(t)) \leq \nn$ for all $t \in [0 ,\infty).$ We denote the equivalence class of a geodesic ray $b$ by $b(\infty).$
 \end{definition}

 Notice that by Proposition 8.2 in the CAT(0) boundary section of \cite{BH1}, for each $b$ representing an element of $\partial X$, and for each $x' \in X$, there is a unique geodesic ray $b'$ starting at $x'$ with $b(\infty)=b'(\infty).$ Now we describe the \emph{cone topology} on this boundary:
 \begin{definition}(Cone topology) \label{def:wConeTopology}
   Fix a base point $\go$ and let $b$ be a geodesic ray starting at $\go$. A neighborhood basis for $b(\infty)$ is given by sets of the form: 
\[ V_{r, \epsilon}(b(\infty)):=\{ c(\infty) \in \partial X| \, c(0)=\go\,\,\text{and }\,d(b(t), c(t))<\epsilon \, \, \text{for all }\, t < r \}.\]

 The topology generated by the above neighborhood basis is called the $\emph{the cone topology}$ and denoted by $\Con$. When we equip $\partial X$ with $\Con$, we denote it $\partial_\infty X$ and call it \emph{the visual boundary.}  Notice that two geodesic rays are close together in the cone topology if they have
representatives starting at the same point which stay close (are at most $\epsilon$ apart) for a long
time (at least $r$). 
 \end{definition}

While the above definition of the visual boundary $\partial_\infty X$ made a reference to a base point $\go$, Proposition 8.8 in the CAT(0) boundaries section of \cite{BH1} shows that the topology of the visual boundary is a base point invariant. It's also worth mentioning that when $X$ is proper, then the space $\overline{X}=X \cup \partial_\infty X$ is compact.




\subsection{Sublinearly Morse geodesics} \label{subsec:SublinearMorseBoundaries}

In this section we review the definition and properties of $\kappa$-Morse geodesic rays needed for this paper. For further details, see \cite{QRT19}. We fix a function 
\[
\kappa : [0,\infty) \to [1,\infty)
\] 
that is monotone increasing, concave, and sublinear, that is
\[
\lim_{t \to \infty} \frac{\kappa(t)} t = 0. 
\]
Note that using concavity, for any $a>1$, we have
\begin{equation} \label{Eq:Concave}
\kappa(a t) \leq a \left( \frac 1a \, \kappa (a t) + \left(1- \frac 1a\right) \kappa(0) \right) 
\leq a \, \kappa(t).
\end{equation}
It is worth noting that the assumption that $\kappa$ is increasing and concave makes certain arguments
cleaner, otherwise they are not really needed. One can always replace any sublinear function
$\kappa$, with another sublinear function $\kappa'$ so that $\kappa(t) \leq  \kappa'(t) \leq  \cc \kappa(t)$ for some constant $\cc$ and $\kappa'$
is monotone increasing and concave.

\begin{notation}
We fix a base point $\go \in X$ once and for all and denote for all $x \in X$, $||x|| := d(\go, x)$. To simplify notation, we often drop $\Norm{\cdot}$. That is, for $x \in X$, we define
\[
\kappa(x) := \kappa(\Norm{x}). \]
\end{notation}

If $x,y$ are within $\kappa(x)$ of each other, then $\kappa(x)$ and $\kappa(y)$ are multiplicatively the same.

\begin{proposition}[{\cite[Lemma 3.2]{QRT19}}] \label{lemma: constants}
For any $\dd_0>0$, there exists $\dd_1,\dd_2>0$ depending only on $\dd_0$ and $\kappa$ so that for $x,y \in X$, we have 
\begin{center}
    $d(x,y) \leq \dd_0 \kappa(x) \implies \dd_1\kappa(x) \leq \kappa(y) \leq \dd_2 \kappa(x).$
\end{center}
\end{proposition}

\begin{definition} ($\kappa$-contracting)
A geodesic ray $b$ is said to be $\kappa$-contracting if there exists a constant $\nn \geq 0$ such that for any ball $B$ centered at $x$ with $B \cap b =\emptyset$, we have $diam(\pi_b(B))< \nn \kappa(x).$
\end{definition}

\begin{definition}[$\kappa$--neighborhood]  \label{Def:Neighborhood} 
For a closed set $Z$ and a constant $\nn$ define the $(\kappa, \nn)$--neighbourhood 
of $Z$ to be 
\[
\calN_\kappa(Z, \nn) = \Big\{ x \in X \ST 
  d_X(x, Z) \leq  \nn \cdot \kappa(||x||)  \Big\}.
\]

\end{definition} 

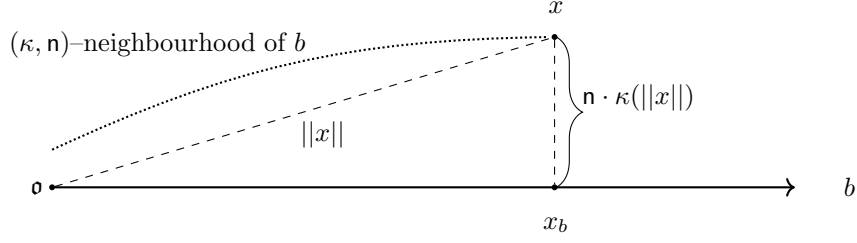
\begin{figure}[h!]
\begin{tikzpicture}
 \tikzstyle{vertex} =[circle,draw,fill=black,thick, inner sep=0pt,minimum size=.5 mm] 
[thick, 
    scale=1,
    vertex/.style={circle,draw,fill=black,thick,
                   inner sep=0pt,minimum size= .5 mm},
                  
      trans/.style={thick,->, shorten >=6pt,shorten <=6pt,>=stealth},
   ]

 \node[vertex] (a) at (0,0) {};
 \node at (-0.2,0) {$\go$};
 \node (b) at (10, 0) {};
 \node at (10.6, 0) {$b$};
 \node (c) at (6.7, 2) {};
 \node[vertex] (d) at (6.68,2) {};
 \node at (6.7, 2.4){$x$};
 \node[vertex] (e) at (6.68,0) {};
 \node at (6.7, -0.5){$x_{b}$};
 \draw [-,dashed](d)--(e);
 \draw [-,dashed](a)--(d);
 \draw [decorate,decoration={brace,amplitude=10pt},xshift=0pt,yshift=0pt]
  (6.7,2) -- (6.7,0)  node [black,midway,xshift=0pt,yshift=0pt] {};

 \node at (7.8, 1.2){$\nn \cdot \kappa(||x||)$};
 \node at (3.6, 0.7){$||x||$};
 \draw [thick, ->](a)--(b);
 \path[thick, densely dotted](0,0.5) edge [bend left=12] (c);
\node at (1.4, 1.9){$(\kappa, \nn)$--neighbourhood of $b$};
\end{tikzpicture}
\caption{A $\kappa$-neighbourhood of a geodesic ray $b$ with multiplicative constant $\nn$.}
\end{figure}

\begin{definition}[$\kappa$--fellow travelling] \label{Def:Fellow-Travel}
Let $\alpha$ and $\beta$ be two quasi-geodesic rays in $X$. If $\alpha$ is contained in some 
$\kappa$--neighbourhood of $\beta$ and $\beta$ is contained in some 
$\kappa$--neighbourhood of $\alpha$, we say that $\alpha$ and $\beta$ 
\emph{$\kappa$--fellow travel} each other. This defines an equivalence
relation on the set of quasi-geodesic rays in $X$.
\end{definition}

\begin{definition}[$\kappa$--Morse] \label{Def:Morse} 
A set $Z$ is \emph{$\kappa$--Morse} if there is a function
\[
\mm_Z : \RR_+^2 \to \RR_+
\]
so that if $\beta : [s,t] \to X$ is a $(\qq, \sQ)$--quasi-geodesic with end points 
on $Z$ then
\[
\beta[s,t]  \subset \calN_{\kappa} \big(Z,  \mm_b(\qq, \sQ)\big). 
\]
We refer to $\mm_{Z}$ as the \emph{Morse gauge} for $Z$. For convenience, we may always assume that $\mm_Z(\qq,\sQ)$ is the largest element in the set $\{\qq,\sQ, \mm_Z(\qq,\sQ) \}$. It is worth noting that when $\kappa=1$ we recover the standard definition of a \emph{Morse set}. For a good reference on Morse geodesics, see \cite{Cordes2017ASO} by Cordes.
\end{definition}

We remark that we will be mostly interested in the case where the set $Z$ is a geodesic ray $b$.

\begin{definition}[Sublinearly Morse boundary] \label{def: SM}
Let $\kappa$ be a sublinear function and let $X$ be a $\CAT$ space. We define the $\kappa$-Morse boundary, as a set, by
\[\partial_\kappa X : = \{ \text{ all } \kappa\text{-Morse quasi-geodesics } \} / \kappa\text{-fellow travelling}\] The $\kappa$-Morse boundary is equipped with the following topology. Fix a base point $\go$, let $\xi \in \partial_{\kappa} X$, and let $b$ be the unique geodesic representative of $\xi$ that starts at $\go$. For all $r > 0$, we define $U_{\kappa}(b, r)$ to be the set of all points $\eta \in \partial_{\kappa} X$ such that for every $(\qq,\sQ)$-quasi-geodesic $\beta$ representing $\eta$, starting at $\go$, and satisfying $m_{b}(\qq,\sQ) \leq \frac{r}{2\kappa(r)}$, we have $\beta\vert_{[0, r]} \subset \mathcal{N}_{\kappa}(b, m_{b}(\qq, \sQ))$ (see \cite{QRT19}). We will denote the topology induced by these neighbourhood bases by $\SM$.

\end{definition}

We compare the $\SM$-topology to the cone topology induced by the visual boundary:

\begin{lemma}\label{lem:injection} ($\Con \subseteq \SM$ over $\kappa$-Morse rays)
There is an injective map $i$ of the set of points in $\partial_\kappa X$ into the visual boundary $\partial_\infty X$. Furthermore, the map $i: \partial_\kappa X \rightarrow \partial_\infty X$ is continuous. When equipping the image of this map with the subspace topology we denote it $\partial^\kappa _\infty X$. In other words, $\partial^\kappa _\infty X$ is the subspace of $\partial_\infty X$ consisting of all $\kappa$-Morse geodesic rays emanating from a fixed based point. 
\end{lemma}

 We note that these topologies are likely not equal in general. This stems from the fact that $\SM$ has a lot of similarities with the topology $\FQ$ that Cashen-Mackay introduced for the Morse boundary (see \cite{CashenMackay}). Thus, the example that shows that $\Con \neq \FQ$ in \cite{IM20a} likely generalises to sublinearly Morse boundaries.

\begin{proof} The proof of the existence of the set map is immediate by the first part of Theorem \ref{thm:Morse iff contracting iff excursion}.
Note that, if two topologies $\mathcal{T}$, $\mathcal{T}'$ are defined by neighbourhood bases, we have that, whenever every basic neighbourhood $U$ of a point $x$ with respect to $\mathcal{T}$ admits a basic neighbourhood $U'$ of $x$ with respect to $\mathcal{T}'$ such that $U' \subseteq U$, then $\mathcal{T}' \supseteq \mathcal{T}$. Our strategy is to find a suitable neighbourhood basis for the topology $\Con$.

Fix a base point $\go \in X$. Let $\xi \in \partial_{\kappa} X$, $R, \epsilon > 0$, and $\xi_{\go}$ the unique geodesic representative of $\xi$ starting at $\go$. By definition of the cone topology, the sets
\[ V_{R, \epsilon}(\xi) := \{ \eta \in \partial_{\kappa}X \vert d( \eta_{\go}(R), \xi_{\go}(R) ) < \epsilon \} \]
define a neighbourhood basis of the cone topology pulled back via the inclusion $i$. We also define the set $V_{\kappa}(\xi_{\go}, r)$ to be the set of all points $\eta \in \partial_{\kappa} X$ such that $\eta_{\go} \vert_{[0, r]} \subset \mathcal{N}_{\kappa}(\xi_{\go}, m_{\xi_{\go}}(1, 0))$. Clearly, for all $r$ satisfying $m_{\xi_{\go}}(1, 0) \leq \frac{r}{2\kappa(r)}$, we have $U_{\kappa}(\xi_{\go}, r) \subseteq V_{\kappa}(\xi_{\go}, r)$. We conclude that, if the sets $V_{\kappa}(\xi_{\go}, r)$ form a neighbourhood basis for $\Con$, then $\FQ \supseteq \Con$.

We are left to show that, for $r$ large, we find $R'$, such that $V_{\kappa}(\xi_{\go}, r) \supseteq V_{R', 1}(\xi)$, as this would prove that the $V_{\kappa}(\xi_{\go}, r)$ are a neighbourhood basis for $\Con$. (The other inclusion is immediate as $\kappa \geq 1$.) Let $\beta$ be a geodesic ray starting at $\go$, respresenting a point $\eta \in \partial_{\kappa} X$ and suppose that $\beta\vert_{[0, r]} \nsubseteq \mathcal{N}_{\kappa}(\xi_{\go}, m_{\xi_{\go}}(1, 0))$. Since $\kappa$-neighbourhoods in $\CAT$ spaces are convex, this implies that
\[ d(\beta(r), \xi_{\go}(r)) > m_{\xi_{\go}}(1, 0) \kappa( \Vert \beta(r) \Vert). \]
Using convexity of distance functions, we conclude that for all $r' \geq r$, we have
\[ d(\beta(r'), \xi_{\go}(r') ) > \frac{r'}{r} \cdot m_{\xi_{\go}}(1, 0) \kappa( \Vert \beta(r) \Vert). \]
Choosing $r'$ sufficiently large, we find $R'$ such that 
\[ \frac{R'}{r} \cdot m_{\xi_{\go}}(1, 0) \kappa( \Vert \beta(r) \Vert) > 1. \]
We conclude that for all $\eta \notin V_{\kappa}(\xi_{\go}, r)$, we have $\eta \notin V_{R', 1}(\xi)$. Therefore, $V_{R', 1}(\xi) \subseteq V_{\kappa}(\xi_{\go}, r)$. This implies that the sets $V_{\kappa}(\xi_{\go}, r)$ form a neighbourhood basis of $\Con$, and $\SM \supseteq \Con$.
\end{proof}







\subsection{CAT(0) cube complexes} \label{subsec:CCCs}

For a more detailed introduction to $\CAT$ cube complexes, see \cite{Sageev}. Our goal for this section is to fix notation and to recall some definitions and facts that will be used in the remainder of the paper.

\begin{definition}(Cubes and midcubes)  A \emph{cube} is a Euclidean unit cube 
$[0,1]^n$
for some
$n \geq  0$. A \emph{midcube} of c is a subspace obtained by restricting exactly one coordinate in $[0,1]^n$ to $\frac{1}{2}$.

\end{definition}

\begin{definition}[Cube complexes]
For $n \geq 0$, let $[0,1]^n$ be an $n$-cube equipped with the Euclidean metric. We obtain a \emph{face} of a $n$--cube by choosing some indices in $\{ 1, \dots, n \}$ and considering the subset of all points where we for each chosen index $i$, we fix the $i$-th coordinate either to be zero or to be one. A \emph{cube complex} is a topological space obtained by glueing cubes together along faces, i.e. every gluing map is an isometry between faces. One may think of a cube complex as a CW-complex whose cells are equipped with the geometry of Euclidean cubes and whose glueing maps are isometries of faces.
\end{definition}

\begin{definition} (local/uniform local finiteness) Let $X$ be a $\CAT$ cube complex and $x \in X^{(0)}$ a vertex. The \emph{degree}, denoted $\deg(x)$ of $x$, is defined to be the number of edges incident to $x$. We have the following:
\begin{itemize}
    \item $X$ is said to be \emph{finite dimensional} if there is an integer $v_X$ such that whenever $[0,1]^n$ is a cube in $X$, we have $n \leq v_X.$
    \item $X$ is said to be \emph{locally finite} if every vertex $x \in X^{(0)}$ has finite degree.
    \item $X$ is said to be \emph{uniformly locally finite} if there is an integer $u_X$ such that for every vertex $x \in X^{(0)}$, we have $\deg(x) \leq u_X$.
\end{itemize}

\end{definition}

Any cube complex can be equipped with a metric as follows. The $n$--cubes are equipped with the Euclidean metric, which allows us to define the length of continuous paths inside the cube complex. (Simply partition every path into segments that lie entirely within one cube and use the Euclidean metric of that cube.) We define
\[ d(x,y) := \inf \{ length(\gamma) \vert \gamma \text{ a continous path from $x$ to $y$} \}. \]
The map $d$ defines a metric on $X$. We sometimes call $d$ the metric induced by the Euclidean metric on each cube.

\begin{definition}[$\CAT$ cube complexes]
Let $X$ be a cube complex and $d$ the metric induced by the Euclidean metric on each cube. We say that $X$ is a $\CAT$ cube complex if and only if $(X,d)$ is a $\CAT$ space.

\end{definition}

\begin{definition}(Hyperplanes, half spaces and separation) Let $X$ be a $\CAT$ cube complex. A \emph{hyperplane} is a connected
subspace $h \subset X$ such that for each cube $c$ of $X$, the intersection $h \cap c$ is either empty or a
midcube of $c$.  For each hyperplane $h$, the complement $X \setminus h$ has exactly two components, called \emph{halfspaces}
associated to $h$. We usually denote these $h^+$ and $h^-$ . If $h^+$ is a halfspace, there is a unique
hyperplane $h$ such that $h^+$ is a component of $X \setminus h$, and we also say $h$ is the hyperplane
associated to $h^+$ . A hyperplane $h$ is said to \emph{separate} the sets $U,V \subseteq X$ if $U \subseteq h^+$ and $V \subseteq h^-.$
\end{definition}

Whenever $X$ is a $\CAT$ cube complex, we refer to the metric induced by the Euclidean metric on each cube as the \emph{CAT(0)-metric}. In contrast to the $\CAT$ metric, we can also equip every $n$--cube with the restriction of the $l^1$ metric of $\mathbb{R}^n$ and consider the induced path metric $d^{(1)}(\cdot, \cdot)$. We refer to $d^{(1)}$ as the {\it combinatorial metric} (or {\it $l^1$--metric}). The following lemma is a standard lemma about how the $l^1$-metric relates to the CAT(0) metric of a CAT(0) cube complex, for example, see \cite{CapriceSageev}.

\begin{lemma} \label{lem: quasi-isometry between combinatorial and CAT(0) metrics}
If $X$ is a finite-dimensional $\CAT$ cube complex, then the $\CAT$ metric $d$ and the combinatorial metric $d^{(1)}$ are bi-Lipschitz equivalent and complete. In particular, if all cubes in $X$ have dimension $\leq m$, then $d \leq d^{(1)} \leq \sqrt{m} d$. Furthermore, for two vertices $x,y \in X^{(0)}$, we have $$d^{(1)}(x,y)=|\{\text{hyperplanes $h \subseteq X$ which separate the vertices }x,y\}|.$$
\end{lemma}

Since we now have two metrics on the $\CAT$ cube complex $X$, we need to distinguish between terminologies for the two metrics. A path in $X$ is called a geodesic when it is a geodesic with respect to the $\CAT$ metric.

\begin{definition}(Combinatorial geodesics)
A path in the $1$--skeleton of $X$ is called a \emph{combinatorial geodesics} if it is a geodesic between vertices of $X$ with respect to the combinatorial metric. {\it Combinatorial geodesic rays} are defined analogously.
\end{definition}

If all cubes in $X$ have dimension at most $m$, every combinatorial geodesic is a $(\sqrt{m}, 0)$-quasi-geodesic. Combinatorial geodesics are fully determined by their starting vertex and the ordered sequence of hyperplanes they cross.\\

\begin{lemma}[Proposition 2.8 in \cite{BeyFio}] \label{lem:ApproximationofGeodesicRays}
Let $X$ be a finite-dimensional $\CAT$ cube complex of dimension at most $m$. Every $\CAT$ geodesic ray $b$ starting at a vertex $\go$ is at Hausdorff-distance at most $m$ from a combinatorial geodesic $\alpha$ with the same origin. Furthermore, $\alpha$ can be chosen such that it crosses the same hyperplanes as $b$.
\end{lemma}

\begin{notation}
We adopt the following notations:

\begin{itemize}
    \item We denote the set of all hyperplanes and of all halfspaces of $X$ by $\mathcal{W}(X)$ and $\mathcal{H}(X)$ respectively.
    \item Given a halfspace $s$, denote the other halfspace bounded by the same hyperplane by $s^*$.
    \item Given a hyperplane $h$, we call a choice of halfspace bounded by $h$ an \emph{orientation of $h$}.
    \item We sometimes denote the two orientations of $h$ by $\{ h^{+}, h^{-} \}$.
    \item Given two subsets $U$, $V \subset X$, we define $\mathcal{W}(U \vert V)$ to be the set of all hyperplanes that separate $U$ from $V$. Given a (combinatorial) geodesic $\alpha$, we write $\mathcal{W}(\alpha)$ for the set of hyperplanes crossed by $\alpha$.
    
\end{itemize}

\end{notation}

When $X$ is fixed, we will simply use $\mathcal{H}$ instead of $\mathcal{H}(X).$ The set $\mathcal{H}$ is endowed with the order relation given by inclusions; the involution $*$ is order reversing. The triple $(\mathcal{H},\cu,*)$ is thus a \emph{pocset} (see \cite{Sageev}).

\begin{definition}\label{def: chain} (Chain). A \emph{chain} in $X$ is a (possibly infinite) collection of mutually disjoint hyperplanes which are associated to a collection of nested half-spaces. 

\end{definition}

\begin{definition} (transverse hyperplanes)
Two hyperplanes are called \emph{transverse} if they intersect. Analogously, two halfspaces $s$, $s' \in \mathcal{H}$ are called \emph{transverse}, if their bounding hyperplanes are transverse. Equivalently, they are transverse if and only if the four intersections $s \cap s'$, $s \cap s'^*$, $s^* \cap s'$, $s^* \cap s'^*$ are non-empty. A hyperplane is \emph{transverse} to a halfspace $s$ if it is transverse to the hyperplane that bounds $s$.
\end{definition}

Note that every intersection $h_1\cap \dots \cap h_k$ of pairwise transverse hyperplanes $h_1, \dots, h_k$ inherits a $\CAT$ cube complex structure. Its cells are precisely the intersections $h_1\cap \dots \cap h_k\cap c$ for any cube $c\cu X$. Alternatively, $h_1\cap \dots \cap h_k$ can be viewed as a subcomplex of the cubical subdivision of $X$.


\begin{definition} (Convex subcomplexes) Let $X$ be a CAT(0) cube complex. We define the following:

\begin{enumerate}
    \item  A subcomplex $Y \subseteq X$ is said to be \emph{a convex subcomplex} with respect to the CAT(0) metric if whenever $x,y \in Y$, we have $[x,y] \in Y$ where $[x,y]$ is the unique CAT(0) geodesic connecting $x,y.$
    \item  A subcomplex $Y \subseteq X$ is said to be \emph{a convex subcomplex} with respect to the combinatorial metric if whenever $x,y \in Y$, we have $[x,y]_1 \in Y$ where $[x,y]_1$ denotes the collection of combinatorial geodesics from $x$ to $y.$
\end{enumerate}

\end{definition}

\begin{definition}(Convex hulls) \label{def: convex hull} Let $\alpha$ be a combinatorial geodesic ray in a $\CAT$ cube complex $X.$ The \emph{convex hull} of $\alpha,$ denoted Hull$(\alpha)$, is defined to be the intersection of all half spaces which contain the combintorial geodesic ray $\alpha$. Notice that for any combinatorial geodesic ray $\alpha$, the space $Hull(\alpha)$ is convex in both the $\CAT$ and the combinatorial metric. An essential property of a convex hull of a combinatorial geodesic ray $\alpha$ is that a hyperplane $h$ meets $Hull(\alpha)$ if and only if $h$ meets $\alpha.$

\end{definition}

The following is a standard lemma about CAT(0) cube complexes.

\begin{lemma}(the Helly property) Any finite collection of convex subcomplexes of a finite dimensional CAT(0) cube complex $X$ satisfy the Helly property. That is to say, any finite collection of pairwise intersecting convex subcomplexes has nonempty total intersection.
\end{lemma}

\begin{definition}(Combinatorial projection) Given a subset $S \subset X$ that is closed and convex with respect to the combinatorial metric, we can define the \emph{combinatorial projection} $\sg_S : X \rightarrow S$, which sends every point $x \in X$ to the unique point in $S$ that is closest to $S$. If $x$ is a vertex and $S$ is a convex subcomplex, then $\sg_S(x)$ can be characterised as the unique vertex $s \in S$, such that for every hyperplane $h \in \mathcal{W}(S)$ the points $s, x$ are contained in the same halfspace bounded by $h$. That is to say, a hyperplane $h$ separates $x$ from it's combinatorial $\sg_S(x)$ if and only if $h$ separates $x$ from $S$. Since hyperplanes are closed and convex, they in particular have a combinatorial projection

\end{definition}

.\\

\begin{definition}($\Hyp$-topology on $\partial X$) For a CAT(0) cube complex $X$, recall that $\partial X$ denotes the set of all equivalence classes of CAT(0) geodesic rays modulo finite Hausdorff. We define the following topology on the set $\partial X$: Fix a vertex $\go \in X$ as a base point and let $h_1, \dots, h_n$ be distinct hyperplanes in $X$. Define

\begin{equation*}
\begin{split}
V_{\go, h_1, \dots, h_n} := \{ \xi \in \partial X | \text{ The unique} & \text{ geodesic representative of $\xi$ based}\\
& \text{ at $\go$ crosses the hyperplanes $h_1, \dots, h_n$} \}.
\end{split}
\end{equation*}

The collection $B=\{ V_{\go, h_1, \dots, h_n}|  n \in \mathbb{N}, \,h_1,h_2,..,h_n \text{ are hyperplanes} \}$ forms the basis of a topology which we denote $\Hyp$.
\end{definition}


We now introduce some notions that are helpful to detect Morse and $\kappa$-Morse geodesic rays in $\CAT$ cube complexes. Two hyperplanes are said to be \emph{$k$-separated} if the number of hyperplanes crossing both of them is bounded from above by $k$. In \cite{ChSu2014}, Charney and Sultan give the following characterization for $1$--Morse geodesic rays in $\CAT$ cube complexes.

\begin{theorem}[{\cite[Theorem 4.2]{ChSu2014}}]\label{thm: Charney and Sultan hyperplanes crossing}
Let $X$ be a uniformly locally finite \CAT \text{ }cube complex and $c \geq 0$. There exist
$r > 0,$ $k \geq 0$ (depending only on $c$ and $X$) such that a geodesic ray $b$ in $X$ is $1$--Morse with Morse gauge $c$ (cf. Definition \ref{Def:Morse})
if and only if $b$ crosses an infinite sequence of hyperplanes $h_1, h_2, h_3$, . . . at points $x_i =
b \cap h_i$ satisfying:

\begin{enumerate}
    \item  $h_i
, h_{i+1}$ are $k$-separated and
\item $d(x_i,x_{i+1}) \leq r.$
\end{enumerate}

\end{theorem}

Since $1$--Morse geodesic rays in $\CAT$ spaces always admit a constant function as a Morse gauge, this provides a characterisation of all $1$--Morse geodesics in $\CAT$ cube complexes. We now provide an analogous characterisation of $\kappa$-Morse geodesic rays.

\begin{definition}[Facing Triples]\label{def:facing triples}
A collection of three hyperplanes $h_1, h_2, h_3$ is said to form a \emph{facing triple} if they are disjoint and none of the three hyperplanes separates the other two.
\end{definition}

Notice that if a (combinatorial or $\CAT$) geodesic $b$ crosses three disjoint hyperplanes $h_1,h_2$, and $h_3$ in that order, then $h_2$ separates $h_1$ and $h_3.$ In particular, a geodesic $b$ cannot cross a facing triple. Conversely, if $\mathcal{C}$ is a collection of hyperplanes  which contains no facing triple, then there is a geodesic which crosses a \emph{definite proportion} of the hyperplanes in $\mathcal{C}$:
    
\begin{lemma}[{\cite[Corollary 3.4]{Hagen20}}]\label{lem:definite proportion} Let $X$ be a $\CAT$ cube complex of dimension at most $m$. There exists a constant $k$, depending only on $m$, such that the following holds.  If $\mathcal{C}$ is a collection of hyperplanes which contains no facing triple, then there exists a (combinatorial or $\CAT$) geodesic which crosses at least $\frac{\mathcal{|C|}}{k}$ hyperplanes from the collection $\mathcal{C}.$
\end{lemma}
 
\begin{definition}[Well-separated hyperplanes]
Two disjoint hyperplanes $h_1, h_2$ are said to be \emph{$k$-well-separated} if any collection of hyperplanes intersecting both $h_1, h_2$, which does not contain any facing triple, has cardinality at most $k$. We say that $h_1$ and $h_2$ are \emph{well-separated} if they are $k$-well-separated for some $k$.
\end{definition}

\begin{theorem}[\cite{Murray-Qing-Zalloum}, Theorem B] \label{thm: characterizing sublinear Morse geodesics using hyperplanes}
Let $X$ be a finite dimensional $\CAT$ cube complex. A geodesic ray $b \in X$ is $\kappa$-contracting if and only if there exists $\cc>0$ such that $b$ crosses an infinite sequence of hyperplanes $h_1, h_2,...$ at points $b(t_i)$ satisfying: 

\begin{enumerate}
    \item $d(t_i,t_{i+1}) \leq \cc \kappa(t_{i+1}).$
    \item $h_i,h_{i+1}$ are $\cc \kappa(t_{i+1})$-well-separated.
\end{enumerate}

\end{theorem}

\begin{definition} [$\kappa$-excursion geodesics]
A geodesic ray $b \in X$ is said to be a \emph{$\kappa$-excursion geodesic} if there exists $\cc>0$ such that $b$ crosses an infinite sequence of hyperplanes $h_1, h_2,...$ at points $b(t_i)$ satisfying: 

\begin{enumerate}
    \item $d(t_i,t_{i+1}) \leq \cc \kappa(t_{i+1}).$
    \item $h_i,h_{i+1}$ are $\cc \kappa(t_{i+1})$-well-separated.
\end{enumerate}
The constant $\cc$ will be referred to as the \emph{excursion constant} and the hyperplanes $\{h_i\}$ will be referred to as the \emph{excursion hyperplanes}.
\end{definition}

\begin{theorem} \label{thm:Morse iff contracting iff excursion} 
Let $X$ be a proper $\CAT$ space and let $\zeta$ be an equivalence class of $\kappa$-fellow travelling quasi-geodesics in $X$. If $\zeta$ contains a $\kappa$-Morse quasi-geodesic ray, then $\zeta$ contains a unique geodesic ray $b$ starting at $b(0)$. Furthermore, if $\eta$ is an equivalence class of $\kappa$-fellow travelling quasi-geodesics, then the following are equivalent:

\begin{itemize}
    \item $\eta$ contains a $\kappa$-Morse geodesic ray.
    \item $\eta$ contains a $\kappa$-contracting geodesic ray.
    \item  $\eta$ contains a $\kappa$-excursion geodesic ray.

    \item Every quasi-geodesic ray in $\eta$ is $\kappa$-Morse.
     \item There exists a quasi-geodesic ray in $\eta$ which is $\kappa$-Morse.

\end{itemize}
In particular, a geodesic ray $b$ is $\kappa$-Morse if and only if it is $\kappa$-contracting if and only if it is a $\kappa$-excursion geodesic.

\end{theorem}

\begin{proof}
The statement that if $\zeta$ contains a $\kappa$-Morse quasi-geodesic ray, then $\zeta$ contains a unique geodesic ray $b$ starting at $b(0)$ is precisely Lemma 3.5 and Proposition 3.10  of \cite{QRT19}. The rest of the statement is exactly Theorem 3.8 in \cite{QRT19} and Theorem \ref{thm: characterizing sublinear Morse geodesics using hyperplanes}.

\end{proof}

\subsection{The Roller boundary} \label{subsec:RollerBoundary}

Recall that for a CAT(0) cube complex $X$, the set $\mathcal{H}(X)$, or simply $\mathcal{H}$, is the set of all half spaces of $X$. Given a vertex $x\in X^{(0)}$, we denote by $\sigma_x \cu \mathcal{H}$ the set of all halfspaces containing the point $x$.

\begin{definition}(The Roller boundary) Let $\iota : X \rightarrow \prod_{h \in \mathcal{W}} \{ h^{+}, h^{-} \}$ denote the map that takes each point $x$ to the set of all half spaces containing $x$ given by $\sigma_x$. Endowing $\prod_{h \in \mathcal{W}} \{ h^{+}, h^{-} \}$ with the product topology, we can consider the closure $\overline{\iota(X)}$. Equipped with the subspace topology, the space $\overline{\iota(X)}$ is a compact Hausdorff space known as the \emph{Roller compactification} of $X$. The \emph{Roller boundary} $\partial_R X$ is defined as the difference $\overline {\iota(X)}\setminus \iota(X)$. If $X$ is locally finite, $\iota(X)$ is open in $\overline X$ and $\partial_R X$ is compact; however, this is not true in general.

\end{definition}

The idea of the following definition is to abstract the properties met by the collection of hyperplanes $\sigma_x$ for some $x \in X^{(0)}$

\begin{definition}(Ultrafilters) A collection of half spaces $\sigma \subseteq \mathcal {H}$ is called a \emph{DCC ultrafilter} if it satisfies the following three conditions:

\begin{enumerate}
\item given any two halfspaces $s$, $s' \in \sigma$, we have $s \cap s' \neq \emptyset$;
\item for any hyperplane $h \in \mathcal{W}$, a side of $h$ lies in $\sigma$;
\item every descending chain of halfspaces in $\sigma$ is finite. 
\end{enumerate}
We refer to a set $\sigma \cu \mathcal{H}$ satisfying only $(1)$ and $(2)$ simply as an \emph{ultrafilter}.

\end{definition}

\begin{remark}
   We remark that the image  $\iota (X)$ is precisely the collection of all DCC ultrafilters. Furthermore, $\overline{\iota(X)}$ coincides with the set of all ultrafilters. We prefer to imagine $\partial_R X$ as a set of points at infinity, represented by combinatorial geodesic rays in $X$, rather than a set of ultrafilters. We will therefore write $x \in \partial_R X$ for points in the Roller boundary and employ the notation $\sigma_x \cu \mathcal{H}$ to refer to the ultrafilter representing $x$.\\
\end{remark}

We say that a hyperplane $h$ \emph{separates} two points $x$ and $y$ in the Roller boundary if $\sigma_x$ and $\sigma_y$ do not contain the same halfspace bounded by $h$. In other words, they induce opposite orientation on $h$. Recall that for subsets $U,V$ of $X$, we use the notation $\mathcal{W}(U \vert V)$ to denote the set of all hyperplanes that separate $U$ from $V$.

\begin{definition}[Combinatorial Gromov product, cf. \cite{BeyFioMedici}] \label{def:CombinatorialGromovProduct}
Let $x, y \in \partial_R X$ and $\go \in X^{(0)}$ be a vertex. We define the \emph{combinatorial Gromov product} of $x$ and $y$ with respect to $\go$ to be
\[ [ x \vert y ]_\go := \# \mathcal{W}(\go \vert x, y), \]
i.e.\,the number of hyperplanes $h$, for which $\sigma_\go$ does not contain the same halfspace bounded by $h$ as $\sigma_x$ and $\sigma_y$.
\end{definition}

\begin{definition}(Components) \label{def:RollerComponents}
Two points $x$, $y \in \partial_R X$ lie in the same \emph{component} if and only if there are only finitely many hyperplanes that separate $x$ from $y$. This defines an equivalence relation on $\partial_R X$ and partitions the Roller boundary into equivalence classes, called components. We denote the component containing a point $x \in \partial_R X$ by $C(x)$. 
\end{definition}

Each component inherits the structure of a $\CAT$ cube complex whose hyperplanes are a strict subset of the set of hyperplanes of $X$. We say that a hyperplane $k \in \mathcal{W}(X)$ \emph{intersects} a component $C$ whenever it corresponds to a hyperplane in $C$. Note that for any two hyperplanes $h$, $k$ that intersect a component $C$, there exist infinitely many $h_i \in \mathcal{W}(X)$ that intersect both $h$ and $k$.

Given two components $C(x), C(y)$, we say that $C(x) \prec C(y)$ if and only if the combinatorial Gromov product $[x \vert y]_\go = \infty$ and any descending chain of half spaces containing $y$ and not containing $x$ terminates. (This partial order was introduced by Guralnik. See for example the appendix of \cite{Genevois2020} for more information.)

\begin{definition}(Median) Let $x, y, z \in X^{(0)} \cup \partial_R X$. We define the \emph{median} of the triple $x, y, z$ to be the unique point $m(x,y,z) \in X^{(0)} \cup \partial_R X$ obtained by orienting every hyperplane of $X$ towards its halfspace that contains the majority of the points $x, y, z$. This orientation yields an ultrafilter and therefore a point in the Roller compactification of $X$, making the median well-defined.
\end{definition}

We know consider a special subspace of the Roller boundary:

\begin{definition} ($\kappa$-Morse points in the Roller boundary) \label{def: kappa-Morse roller}
Let $X$ be finite dimensional. A point $x \in \partial_R X$ is called \emph{$\kappa$--Morse} if it admits a combinatorial geodesic representative that is $\kappa$--Morse. Denote the set of all $\kappa$--Morse points in the Roller boundary of $X$ by $\partial_{R}^{\kappa} X$. 

\end{definition}
Since we often deal with different topologies on the same boundary of the CAT(0) space $X$, we make explicit which notations we use for which boundary:

\begin{notation}(Various boundaries) \label{notation: boundaries}

\begin{itemize}
    \item We denote the equivalence classes of the CAT(0) geodesic rays in $X$ as a \emph{set} by $\partial X$.
    \item When equipping the set $\partial X$ with the \emph{cone topology} given in Definition \ref{def: visual topology}, we denote it $\partial_\infty X$ and refer to it as the \emph{visual boundary}.
    
    \item As a subset of the visual boundary, the collection of all $\kappa$-Morse geodesic rays as a \emph{set} is denoted by $\partial^\kappa X$.
    \item When equipping the subset $\partial^\kappa X$ with the \emph{subspace topology} from the visual boundary, we denote it by $\partial^\kappa_\infty X$.
    \item The set $\partial_\kappa X$ is the set of all equivalence classes of $\kappa$-Morse quasi geodesic rays with respect to the $\kappa$-fellow travelling relation from Definition \ref{Def:Fellow-Travel}.
     \item When we equip the set $\partial_\kappa X$ with the $\SM$-topology from Definition \ref{def: SM}, we will still denote it by $\partial_\kappa X.$
     \item The subspace $\partial_R ^\kappa X$ denotes the subspace of the Roller boundary consisting of points that can be represented by $\kappa$-Morse combinatorial geodesic rays (see Definition \ref{def: kappa-Morse roller}).
     
     \item For a geodesic ray $b$, we will often confuse $b$ with its image $im(b)$. Since we deal with both CAT(0) and combinatorial rays, we will use the letters $a,b,c,\dots$ denote CAT(0) geodesic rays, and the letters $\alpha, \beta,\gamma,\dots$ to denote combinatorial geodesic rays.
    
\end{itemize}
\end{notation}







\section{Combinatorial Gromov products and sublinear fellow-traveling} \label{sec:CombinatorialGromovProducts}
In this section, we prove a few useful statements which will be used in the proof of the main theorems. The following lemma is the main lemma of this section, and its essential to prove well-definedness of the map in Theorem \ref{introthm: Roller map}.

\begin{lemma} \label{lem:GromovProductFellowTraveling}
Let $\alpha, \beta$ be two combinatorial geodesic rays in a $\CAT$ cube complex $X$ such that $\alpha(0)=\beta(0) =: \go$ and $[\alpha|\beta]_{\go}=\infty$. If $\alpha$ is $\kappa$-Morse, then $\beta$ is also $\kappa$-Morse. Furthermore, $\alpha$ and $\beta$ $\kappa$-fellow travel.
\end{lemma}

In preparation for the proof of Lemma \ref{lem:GromovProductFellowTraveling}, we prove several other lemmas first.

\begin{lemma} \label{lem: thin squares}
Let $\alpha$ be a combinatorial geodesic ray and let $Y=Hull(\alpha)$. Suppose that $Y$ contains a $\kappa$-excursion CAT(0) geodesic ray $b$. If $\{ h_i\}_{i \in \mathbb{N}}$ is the collection of excursion hyperplanes for $b$, $p_i=h_i \cap im(\alpha)$, and $q_i=im(b) \cap h_i$, then there exist constants $\dd_1,\dd_2$, depending only on $\kappa$, the excursion constant $\cc$ and on the dimension $v_X$ such that:

\begin{enumerate}
    \item $|\mathcal{W}(p_i \vert p_{i+1})|\leq \dd_1 \kappa(t_{i})$,
    \item $|\mathcal{W}(p_i \vert q_i)| \leq \dd_2 \kappa(t_{i})$, 
\end{enumerate}
where $t_{i}=d(q_{i},b(0))$.

\end{lemma}

\begin{proof}
Let $\mathcal{W}_i$ denote the collection of hyperplanes meeting $h_i.$ To simplify notation, denote $\mathcal{W}=\mathcal{W}(p_i \vert p_{i+1}).$ Every hyperplane in $\mathcal{W}$ must live in either $\mathcal{W}(q_i \vert q_{i+1}), \mathcal{W}_i$ or $\mathcal{W}_{i+1}.$ Thus, we have $$|\mathcal{W}|  \leq  |\mathcal{W} \cap \mathcal{W}(q_i \vert q_{i+1})|+ | \mathcal{W} \cap \mathcal{W}_i|+ | \mathcal{W}\cap \mathcal{W}_{i+1}|.$$ Using Lemma \ref{lem: quasi-isometry between combinatorial and CAT(0) metrics},
there exists a constant $\cc'$ such that $|\mathcal{W}(q_i \vert q_{i+1})| \leq \cc'd(q_i,q_{i+1}) \leq \cc'\cc \kappa(t_{i+1})$. Notice that every $h \in \mathcal{W} \cap \mathcal{W}_i$  must live in either $\mathcal{W}_{i-1}$ or $\mathcal{W}(q_{i-1},q_i)$ (see Figure \ref{fig:FivePossibilities}). In other words, we have the following

$$|\mathcal{W} \cap \mathcal{W}_i| \leq |\mathcal{W} \cap \mathcal{W}_i \cap \mathcal{W}_{i-1}|+ |\mathcal{W} \cap \mathcal{W}_i \cap \mathcal{W}(q_{i-1} \vert q_i)|.$$ The collection of hyperplanes $\mathcal{W} \cap \mathcal{W}_i \cap \mathcal{W}_{i-1}$ is one that contains no facing triple since $\mathcal{W} \subseteq W(\alpha)$. Therefore, as $h_{i-1},h_i$ are $\cc \kappa(t_i)$-well-separated, we have $|\mathcal{W} \cap \mathcal{W}_i \cap \mathcal{W}_{i-1}| \leq \cc \kappa(t_i)$. By Lemma \ref{lem: quasi-isometry between combinatorial and CAT(0) metrics}, there exists a constant $\cc'$ such that $|\mathcal{W}(q_{i-1} \vert q_{i})| \leq \cc'd(q_{i-1},q_{i}) \leq \cc'\cc \kappa(t_{i})$. Hence, we have

 \begin{align*}
    |\mathcal{W} \cap \mathcal{W}_i| &\leq |\mathcal{W} \cap \mathcal{W}_i \cap \mathcal{W}_{i-1}|+ |\mathcal{W} \cap \mathcal{W}_i \cap \mathcal{W}(q_{i-1} \vert q_i)|\\
    &\leq \cc\kappa(t_i)+\cc'\cc \kappa(t_i)\\
    &=\cc(1+\cc')\kappa(t_i).
      \end{align*}
      
      Similarly, we have 
       \begin{align*}
       |\mathcal{W} \cap \mathcal{W}_{i+1}| &\leq \cc(1+\cc')\kappa(t_{i+1}).
          \end{align*}
          
          Thus, we have 
          
          \begin{align*}
          |\mathcal{W}|  &\leq  |\mathcal{W} \cap \mathcal{W}(q_i \vert q_{i+1})|+ | \mathcal{W} \cap \mathcal{W}_i|+ | \mathcal{W}\cap \mathcal{W}_{i+1}|\\
          &\leq \cc'\cc\kappa(t_{i+1})+\cc(1+\cc')\kappa(t_i)+\cc(1+\cc')\kappa(t_{i+1})\\
          &\leq \cc'\cc\kappa(t_{i+1})+\cc(1+\cc')\kappa(t_{i+1})+\cc(1+\cc')\kappa(t_{i+1})\\
          &=(\cc'\cc+\cc(1+\cc')+\cc(1+\cc'))\kappa(t_{i+1}).
          \end{align*}
          Letting $\dd_0=\cc'\cc+ 2\cc(1+\cc')$ gives that $|\mathcal{W}(p_i \vert p_{i+1})|\leq \dd_0 \kappa(t_{i+1})$, however, since $t_i, t_{i+1}$ satisfy $d(t_i, t_{i+1}) \leq \cc \kappa(t_{i+1})$, Proposition \ref{lemma: constants} implies that $\kappa(t_{i+1}) \leq \dd_0'\kappa(t_i)$ for a constant $\dd_0'$ depending only on $\kappa$. Taking $\dd_1=\dd_0\dd_0'$ gives that $|\mathcal{W}(p_i \vert p_{i+1})|\leq \dd_1 \kappa(t_{i})$ which finishes the first part of the statement.

          For the second part, let $\mathcal{W}'=\mathcal{W}(p_i \vert q_i)$. Notice that any hyperplane $h \in \mathcal{W}'$ lives in either $\mathcal{W}(q_{i-1} \vert q_{i})$, $\mathcal{W}(p_{i-1} \vert q_{i-1}),$ or $\mathcal{W}(p_{i-1} \vert p_{i})$. In other words, we have

          $$|\mathcal{W}'| \leq |\mathcal{W}' \cap \mathcal{W}(q_{i-1} \vert q_{i})|+|\mathcal{C}' \cap \mathcal{W}(p_{i-1} \vert q_{i-1})|+|\mathcal{C}' \cap \mathcal{W}(p_{i-1} \vert p_{i})|.$$ Now, using part 1 of this lemma and since $b$ is a $\kappa$-excursion geodesic with constant $\cc$ we get:
          
          \begin{align*}
          |\mathcal{W}'| &\leq |\mathcal{W}' \cap \mathcal{W}(q_{i-1} \vert q_{i})|+|\mathcal{W}' \cap \mathcal{W}(p_{i-1} \vert q_{i-1})|+|\mathcal{W}' \cap \mathcal{W}(p_{i-1} \vert p_{i})| \\
          &\leq \cc'\cc \kappa(t_{i+1})+\cc \kappa(t_{i+1})+\dd_0\kappa(t_{i+1})\\
          &= (\cc'\cc +\cc +\dd_0)\kappa(t_{i+1}).
          \end{align*}
          See Figure \ref{fig:ThreePossibilities} for the decomposition of $\mathcal{W}'$. Taking $\cc_0=\cc'\cc +\cc +\dd_0$ gives that   $|\mathcal{W}(p_i \vert q_i)| \leq \cc_0 \kappa(t_{i+1})$. Using Proposition \ref{lemma: constants}, as $d(t_i,t_{i+1}) \leq \cc \kappa (t_{i+1})$, there exists a constant $\cc_0'$, depending only on $\kappa$, such that $\kappa(t_{i+1}) \leq \cc_0' \kappa(t_i) $. Therefore, letting $\dd_2=\cc_0\cc_0'$ gives that $|\mathcal{W}(p_i \vert q_i)| \leq \dd_2 \kappa(t_{i})$.

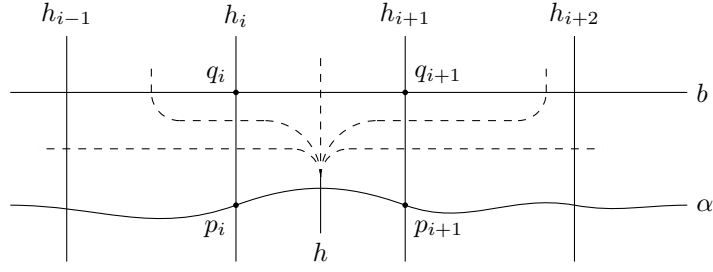
\begin{figure}
\begin{tikzpicture}[scale=1.50]
\draw [] (-3,1) -- (3,1);
\node [right] at (3,1) {$b$};
\draw (-3, 0) to [out = 0, in = 200] (-1, 0);
\draw (-1, 0) to [out = 20, in = 160] (0.5, 0);
\draw (0.5, 0) to [out = -20, in = 170] (2, 0);
\draw (2, 0) to [out = -10, in = 180] (3, 0);
\node [right] at (3,0) {$\alpha$};
\draw (-2.5, -0.5) -- (-2.5, 1.5);
\draw (-1, -0.5) -- (-1, 1.5);
\draw (0.5, -0.5) -- (0.5, 1.5);
\draw (2, -0.5) -- (2, 1.5);
\node [above] at (-2.5, 1.5) {$h_{i-1}$};
\node [above] at (-1, 1.5) {$h_{i}$};
\node [above] at (0.5, 1.5) {$h_{i+1}$};
\node [above] at (2, 1.5) {$h_{i+2}$};
\draw (-0.25, -0.25) -- (-0.25, 0.25);
\node [below] at (-0.25, -0.25) {$h$};
\draw [dashed] (-0.25, 0.25) arc [radius=0.25cm, start angle = 0, end angle = 90];
\draw [dashed] (-0.5, 0.5) -- (-2.75, 0.5);
\draw [dashed] (-0.25, 0.25) arc [radius=0.5cm, start angle = 0, end angle = 90];
\draw [dashed] (-0.75, 0.75) -- (-1.5, 0.75);
\draw [dashed] (-1.5, 0.75) arc [radius=0.25cm, start angle = 270, end angle = 180];
\draw [dashed] (-1.75, 1) -- (-1.75, 1.25);
\draw [dashed] (-0.25, 0.25) -- (-0.25, 1.35);
\draw [dashed] (-0.25, 0.25) arc [radius=0.5cm, start angle = 180, end angle = 90];
\draw [dashed] (0.25, 0.75) -- (1.5, 0.75);
\draw [dashed] (1.5, 0.75) arc [radius=0.25cm, start angle = -90, end angle = 0];
\draw [dashed] (1.75, 1) -- (1.75, 1.25);
\draw [dashed] (-0.25, 0.25) arc [radius=0.25cm, start angle = 180, end angle = 90];
\draw [dashed] (0, 0.5) -- (2.25, 0.5);
\draw [fill] (-1,0) circle [radius=0.02cm];
\node [below left] at (-1,-0.05) {$p_i$};
\draw [fill] (0.5,0) circle [radius=0.02cm];
\node [below right] at (0.5,-0.05) {$p_{i+1}$};
\draw [fill] (-1,1) circle [radius=0.02cm];
\node [above left] at (-1, 1) {$q_i$};
\draw [fill] (0.5,1) circle [radius=0.02cm];
\node [above right] at (0.5, 1) {$q_{i+1}$};
\end{tikzpicture}
\caption{The five possibilities how a hyperplane $h$ separating $p_i$ from $p_{i+1}$ may interact with $h_{i-1}, h_i, h_{i+1}$ and $h_{i+2}$. The number of appearances of each possibility is bounded in terms of $\kappa$ because of $\kappa$-well-separatedness and the fact that the $d(q_j, q_{j+1})$ is $\kappa$-bounded for all $j$.}
\label{fig:FivePossibilities}
\end{figure}

\begin{figure}
\begin{tikzpicture}[scale=1.50]
\draw [] (-1,1) -- (2.5,1);
\node [right] at (2.5,1) {$b$};
\draw (-1, 0) to [out = 20, in = 160] (0.5, 0);
\draw (0.5, 0) to [out = -20, in = 170] (2, 0);
\draw (2, 0) to [out = -10, in = 180] (2.5, 0);
\node [right] at (2.5, 0) {$\alpha$};
\draw (0, -0.5) -- (0, 1.5);
\draw (1.5, -0.5) -- (1.5, 1.5);
\node [above] at (0, 1.5) {$h_{i-1}$};
\node [above] at (1.5, 1.5) {$h_{i}$};
\draw (1.75, 0.5) -- (1.25, 0.5);
\node [right] at (1.75, 0.5) {$h$};
\draw [dashed] (1.25, 0.5) -- (-0.25, 0.5);
\draw [dashed] (1.25, 0.5) arc [radius=0.25cm, start angle = 270, end angle = 180];
\draw [dashed] (1, 0.75) -- (1, 1.25);
\draw [dashed] (1.25, 0.5) arc [radius=0.25cm, start angle = 90, end angle = 180];
\draw [dashed] (1, 0.25) -- (1, -0.25);
\draw [fill] (0, 0.13) circle [radius=0.02cm];
\node [below left] at (0, 0.125) {$p_{i-1}$};
\draw [fill] (1.5,0) circle [radius=0.02cm];
\node [below right] at (1.5,-0.05) {$p_{i}$};
\draw [fill] (0, 1) circle [radius=0.02cm];
\node [above left] at (0, 1) {$q_{i-1}$};
\draw [fill] (1.5,1) circle [radius=0.02cm];
\node [above right] at (1.5, 1) {$q_{i}$};
\end{tikzpicture}
\caption{The three possibilities how a hyperplane $h$ separating $p_i$ from $q_i$ may interact with $p_{i-1}, q_{i-1}$ and $h_{i-1}$. The number of appearances of each possibility is bounded in terms of $\kappa$ because of $\kappa$-well-separatedness, the fact that the $d(q_j, q_{j+1})$ is $\kappa$-bounded for all $j$ and inequality $(1)$.}
\label{fig:ThreePossibilities}
\end{figure}
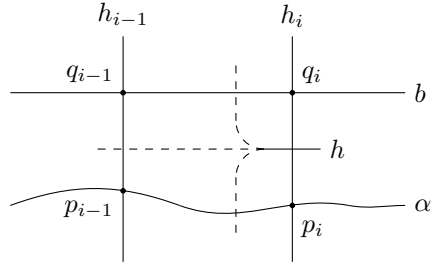

\end{proof}

\newpage

\begin{lemma} \label{lem:BoundofHyperplanesandGeodesics}
Let $b$ be a $\CAT$ $\kappa$-excursion geodesic ray, $\mathcal{H}=\{ h_i\}$ be the collection of excursion hyperplanes crossed by $b$ and let $q_i=im(b) \cap h_i$. There exists a constant $\dd$ such that the following holds:

\begin{enumerate}

\item If $c$ is a geodesic ray crossing $h_{j-1},h_j,h_{j+1} \in \mathcal{H}$ at points $p_{j-1},p_j, p_{j+1},$ then $d(p_j,q_j) \leq \dd \kappa(t_{j}),$ where $t_{j}=d(b(0),q_{j}).$ 
    \item Suppose $h$ is a hyperplane crossing $h_{j-1},h_j,h_{j+1} \in \mathcal{H}$. Let $p_j$ be the closest point projection of $q_j$ to the intersection $h \cap h_j$ and let $p_{j-1}, p_{j+1}$ be points in $h\cap h_{j-1}, h\cap h_{j+1}$ respectively.  We have $d(p_j,q_j) \leq \dd \kappa(t_{j}).$

\end{enumerate}

In particular, if the geodesic $c$ or the hyperplane $h$ crosses infinitely many hyperplanes from $\mathcal{H},$ then $b$ lives in the $\kappa$-neighborhood of $c$ or $h$ respectively.
\end{lemma}

\begin{proof}
The proof of this lemma is very similar to the proof of Lemma \ref{lem: thin squares}, however, for completeness, we give the proof. Let $\cc$ be the excursion constant. Every hyperplane in $W(p_j,q_j)$ must meet either $[q_{j-1},q_j],$ $[p_{j-1}, q_{j-1}]$ or $[ p_{j-1},p_j].$ This implies that 

\begin{equation}\label{E1}
  |W(p_j,q_j)| \leq |W(q_j,q_{j-1})|+|W(q_{j-1},p_{j-1})|+|W(p_{j-1},p_j)|. \tag{1}
\end{equation}

Now we consider the collection $W(p_{j-1},p_j)$. Notice that every hyperplane in $W(p_{j-1},p_j)$ must meet either $[q_j,q_{j+1}]$ or $[q_{j+1}, p_{j+1}]$. That is

\begin{equation}\label{E2}
  |W(p_{j-1},p_j)| \leq |W(q_j,q_{j+1})|+|W(q_{j+1},p_{j+1})|. \tag{2}
\end{equation}

Now, combining the equations \ref{E1} and \ref{E2}, and using the fact that $b$ is a $\kappa$-excursion hyperplane with excursion constant $\cc,$ we get the following:

 \begin{align*}
         |W(p_j,q_j)|  &\leq  |W(q_j,q_{j-1})|+|W(q_{j-1},p_{j-1})|+|W(p_{j-1},p_j)|\\
          &\leq \cc \kappa(t_j)+\cc\kappa(t_j)+|W(p_{j-1},p_j)|\\
          &\leq 2\cc \kappa(t_j)+\cc \kappa(t_{j+1}) +\cc \kappa(t_{j+1})\\
          &=4\cc \kappa(t_{j+1}).
          \end{align*}
Now, using Proposition \ref{lemma: constants}, there exists a constant $\dd_0$ depending only on $\cc$ and $\kappa$ such that $\kappa(t_{j+1}) \leq \dd_0 \kappa(t_j)$. Taking $\dd_1=4\cc \dd_0 $ gives that 
$$|W(p_j,q_j)| \leq \dd_0 \kappa(t_j).$$ However, since the CAT(0) distance is bi-lipschitz equivalent to the $d^{(1)}$-distance (Lemma \ref{lem: quasi-isometry between combinatorial and CAT(0) metrics}), the result follows. 

The proof of part 2 is identical with a hyperplane $h$ replacing the role of the geodesic $c$.

\end{proof}

 \begin{lemma}\label{lem: convexhull contains Morse implies Morse}
Let $\alpha$ be a combinatorial geodesic ray in a $\CAT$ cube complex $X$ and let $Y=Hull(\alpha)$. If $Y$ contains a $\CAT$ geodesic ray $b$ which is $\kappa$-Morse, then the $\CAT$ quasi-geodesic $\alpha$ and $b$ must $\kappa$-fellow travel. In particular, $\alpha$ is $\kappa$-Morse.
\end{lemma}

\begin{proof} Using Theorem \ref{thm:Morse iff contracting iff excursion}, the geodesic ray $b$ must be $\kappa$-excursion. Part 2 of Lemma \ref{lem: thin squares} gives that $\alpha$ and $b$ $\kappa$-fellow travel each other. In other words, both $\alpha$ and $b$ live in the same equivalence class $\sB$. Since $b$ is $\kappa$-Morse, Theorem \ref{thm:Morse iff contracting iff excursion} implies that $\alpha$ is $\kappa$-Morse.

\end{proof}

The following corollary is Corollary 3.5 in \cite{Hagen20}, it's a strengthened version of Lemma 2.1 in \cite{CapriceSageev}.

\begin{corollary} \label{cor: chains in geodesics} (Chains in geodesics) Let $X$ be a $\CAT$ cube complex of dimension $v_X< \infty$. For any $k \in \mathbb{N}$, if $\alpha$ is a
geodesic (in the combinatorial or CAT(0) metric) that crosses at least $v_X \cdot k$ hyperplanes, then
the set of hyperplanes crossing $\alpha$ contains a chain of cardinality $k$. In particular, every geodesic line crosses a bi-infinite chain of hyperplanes.

\end{corollary}

We use the above corollary to show the following.

 \begin{lemma} \label{lem: Existence of a visibility point implies uniqueness}
Let $\alpha$ be a combinatorial geodesic ray in a $\CAT$ cube complex $X$ and let $Y=Hull(\alpha)$. If $\partial_\infty Y$ contains a visibility point, then $|\partial_\infty Y=1|.$
\end{lemma}

\begin{proof}
Recall from Definition \ref{def: convex hull} that a hyperplane $h$ meets $Hull(\alpha)$ if and only if it meets $\alpha$. We first claim that, since $Y = Hull(\alpha)$, there exists no bi-infinite chain $(h_n)_{n \in \mathbb{Z}}$ of hyperplanes in $Y$. For every hyperplane in $Y$, we can choose an orientation (i.e. an associated halfspace) that contains the vertex $\alpha(0)$. Since $\alpha(0) \in X^{(0)}$, we know that every descending sequence of halfspaces in this orientation terminates in finite time. However, since the sequence $(h_n)_{n \in \mathbb{Z}}$ is crossed by the combinatorial geodesic ray $\alpha$, orienting all of them towards $\alpha(0)$ forms (after reordering the hyperplanes appropriately) an infinite descending chain of halfspaces, a contradiction. We conclude that $Y$ contains no bi-infinite chain of hyperplanes.

We now prove the lemma. Let $\zeta$ be the visibility point in  $ \partial_\infty Y$ and let $b$ be a $\CAT$ geodesic ray in $Y$ with $[b]=\zeta$. Suppose for the sake of contradiction that there exists a geodesic ray  $b'$ in $Y$ such that $b'(0)=b(0)$ but $b \neq b'.$ Since $\zeta$ is a visibility point, there exists a geodesic line $l$ connecting $[b]$ to $[b']$. Since $l$ lies in $Y$, every hyperplane meeting $l$ also meets the combinatorial geodesic $\alpha$. Since $l$ is a bi-infinite geodesic line in a finite dimensional $\CAT$ cube complex, using Corollary \ref{cor: chains in geodesics}, the line $l$ must cross a bi-infinite sequence of mutually disjoint hyperplanes. This contradicts the claim that we proved in the first half of this proof and implies the lemma.
\end{proof}

\begin{proposition}[{\cite[Proposition 3.10 ]{QRT19}}]\label{prop: the existence of a fellow travelling geodesic ray} If $\alpha:[0, \infty) \rightarrow X$ is a $\kappa$-Morse quasi geodesic ray, then there exists a geodesic ray $b$ such that $\alpha$ and $b$ $\kappa$-fellow travel each other. Furthermore, $b$ is $\kappa$-Morse.

\end{proposition}

The following corollary states that $\kappa$-Morse geodesic rays define visibility points in the visual boundary.

\begin{corollary}[{\cite[Corollary 1.2 ]{Zalloum20}}]\label{cor: visibility}
Let $Y$ be a proper \CAT space. If $b$ is a $\kappa$-Morse geodesic ray in $Y$, then $b(\infty)$ is a visibility point in $\partial_\infty Y.$
\end{corollary}

\begin{corollary}\label{cor: combinatorial geodesic ray meets same hyperplanes}
Let $\alpha$ be a $\kappa$-Morse combinatorial geodesic ray in a $\CAT$ cube complex and let $Y=Hull(\alpha)$. There exists a $\kappa$-Morse geodesic ray $b \subset Y$ such that $b(0) = \alpha(0)$ and $b$ and $\alpha$ $\kappa$-fellow travel. Furthermore, $b$ is the unique geodesic ray in $Y$ with $b(0)=\alpha(0)$. 
\end{corollary}

\begin{proof}
Since $\alpha$ is a $\kappa$-Morse combinatorial geodesic in the $d^1$-metric of $X$, it must be a $\kappa$-Morse quasi geodesic in the $\CAT$-metric of $X$. Thus, $\alpha$ is a $\kappa$-Morse quasi-geodesic in $Y,$ and hence by Proposition \ref{prop: the existence of a fellow travelling geodesic ray} applied to $\alpha$ and to $Y,$ there exists a $\kappa$-Morse geodesic ray $b \in Y$ such that $\alpha$ and $b$ $\kappa$-fellow travel each other. Using Theorem \ref{cor: visibility}, since $b$ defines a visibility point in the visual boundary $\partial_\infty Y,$ Lemma \ref{lem: Existence of a visibility point implies uniqueness} gives the desired statement. 

\end{proof}

We are now ready to prove the main lemma of this section. We restate it here.

\begin{lemma}
Let $\alpha, \beta$ be two combinatorial geodesic rays in a $\CAT$ cube complex $X$ such that $\alpha(0)=\beta(0) =: \go$ and $[\alpha|\beta]_{\go}=\infty$. If $\alpha$ is $\kappa$-Morse, then $\beta$ is also $\kappa$-Morse. Furthermore, $\alpha$ and $\beta$ $\kappa$-fellow travel.
\end{lemma}

\begin{proof}
Denote $Y_1=Hull(\alpha)$ and $Y_2=Hull(\beta)$. Since $[\alpha|\beta]_{\go}=\infty$, there exists an infinite sequence of hyperplanes $\{h_i\}$ intersecting both $\alpha$ and $\beta$. For any $i,$ since $h_i \cap Y_1, h_i \cap Y_2, Y_1 \cap Y_2$ are all non-empty, there has to exist a point $x_i \in Y_1 \cap Y_2 \cap h_i$ by the Helly property. Notice that since $x_i \in h_i$ and since $X$ is proper, the sequence $x_i$ is unbounded (its not stuck in any compact subset of $X$). Applying Arzelà–Ascoli to the $\CAT$ geodesics $[\go,x_i]$ in $Y_1 \cap Y_2$ yields a $\CAT$ geodesic ray $b \in Y_1 \cap Y_2$ with $b(0)=\go.$ Applying the previous Corollary to $\alpha$ and $Y_1$, we get that $\alpha$ and $b$ $\kappa$-fellow travel and that $b$ is $\kappa$-Morse. Since $b$ is $\kappa$-Morse, the point defined by $b$ in the visual boundary of $Y_2$ is a visibility point, hence by Lemma \ref{lem: Existence of a visibility point implies uniqueness}, we have $|\partial_\infty Y_2|=1$. Now, since $b \in Y_2$ is $\kappa$-Morse, Lemma \ref{lem: convexhull contains Morse implies Morse} implies that $\beta$ is also $\kappa$-Morse, and that $b$ and $\beta$ $\kappa$-fellow travel. Since $\alpha$ and $b$ $\kappa$-fellow travel and $\beta$ and $b$ also $\kappa$ fellow travel, we conclude that $\alpha$ and $\beta$ $\kappa$-fellow travel, in other words, $[\alpha]=[\beta].$

\end{proof}




\section{Comparing HYP to the subspace topology on sublinearly Morse geodesic rays } \label{sec:ComparingTopologies}

In this section, we prove that the topology $\Hyp$ on the sublinear boundary is the same as the subspace topology induced by the cone topology. We begin with the following Lemma.





\begin{lemma} \label{lem:BasisforHYP}
Let $X$ be a finite dimensional CAT(0) cube complex and let $\go \in X$. The open sets $\{U_{\go,k}\}$ form a basis for the topology $\Hyp$ on the set $\partial ^\kappa X.$
\end{lemma}

\begin{proof}
Suppose that $h_1, \dots ,h_n$ is a finite set of hyperplanes crossed by a $\kappa$-Morse geodesic ray $b$. Let $(k_i)_i$ be a family of excursion hyperplanes for $b$. We wish to find some hyperplane $k$ that doesn't intersect $h_1, \dots, h_n$ and is crossed by $b$ after crossing $h_1, \dots ,h_n$. Suppose by contradiction that every hyperplane crossed by $b$ after crossing $h_1, \dots, h_n$ intersects one of the $h_i$. This implies that one of those hyperplanes, we will denote it by $h$, is crossed by an infinite family $(k_{i_j})_j$ of excursion hyperplanes that are crossed by $b$ after crossing $h$. Note that, since all $k_i$ are mutually disjoint and ordered in the order in which they are crossed by $b$, we have that if $h$ intersects $k_i$ and $k_j$ for $i < j$, then $h$ intersects $k_l$ for all $i < l < j$. Therefore, we have that $h$ intersects all elements of the sequence $(k_i)_i$ for $i$ sufficiently large.

Let $q_i := b \cap k_i$ and let $p_i \in h \cap k_i$ be the closest point to $q_i$. Define $t_i := d(b(0), q_i)$. By Lemma \ref{lem:BoundofHyperplanesandGeodesics}, there exists a constant $c'$, such that $d(q_i, p_i) \leq c' \kappa(t_{i})$ for all $i$ sufficiently large. Since $q_i=b(t_i),$ we have $d(b(t_i), h) \leq c \kappa(t_{i})$, which is a sublinear function. However, the function $t \mapsto d(b(t), h)$ is convex. Since all non-negative, sublinear, convex functions are non-increasing functions, $d(b(t), h)$ has to be non-increasing, which contradicts the fact that $b$ intersects $h$. The lemma follows.

\end{proof}

\begin{theorem} \label{thm:VisualequalsHyp}
Let $X$ be a finite dimensional $\CAT$ cube complex. The restrictions of the cone topology and $\Hyp$ to the set of all geodesic rays that can be represented by $\kappa$-Morse geodes rays, denoted $\partial^{\kappa} X$, are equal.
\end{theorem}

\begin{proof}
Let $\xi \in \partial_{\infty} X$. Fix a base point $\go \in X$. Denote the unique geodesic representative of $\xi$ starting at $\go$ by $\xi_\go$. We denote the basic open sets of the cone topology, restricted to $\partial^{\kappa} X$ by
\[ U_{R, \epsilon}(\xi) := \{ \eta \in \partial_{\infty} X \vert d(\xi_\go(R), \eta_\go(R)) < \epsilon \} \cap \partial^{\kappa} X. \]

We first show that the cone topology is finer than $\Hyp$. Let $U_{\go,k}$ be a basic open set of $\Hyp$ and let $\xi \in U_{\go,k}$. We need to find $R, \epsilon > 0$, such that $\xi \in U_{R, \epsilon}(\xi) \subset U_{\go,k}$.

Since $\xi \in U_{\go,k}$, we have that $\xi_\go$ crosses the hyperplane $k$. Therefore, we find some $R > 0$, such that for all $r \geq R$, $d(\xi_\go(r), k) \geq 1$. We claim that $U_{R, \frac{1}{2}}(\xi) \subset U_{\go,k}$. Suppose, $\eta \in U_{R, \frac{1}{2}}(\xi)$. We conclude that the geodesic from $\xi_\go(R)$ to $\eta_\go(R)$ is contained in the ball $B_{\frac{1}{2}}(\xi_\go(R))$ and cannot intersect $k$. This implies that $k$ separates $\go$ and $\eta_\go(R)$. Therefore, $\eta_\go$ crosses $k$ and $\eta \in U_{\go,k}$. This implies that $U_{\go,k}$ is open with respect to the cone topology. Since the sets $\{ U_{\go,k} \}_k$ generate $\Hyp$, we conclude that the cone topology is finer than $\Hyp$ on $\partial^{\kappa} X$.\\

We now prove that $\Hyp$ is finer than the cone topology on $\partial^{\kappa} X$. Let $\xi \in \partial^{\kappa} X$, let $U_{R, \epsilon}(\xi)$ be a basic open set with respect to the cone topology and let $\eta \in U_{R, \epsilon}(\xi)$. Since for every $\xi \in \partial_{\infty} X$, the sets $\{ U_{R, \epsilon}(\xi) \}_{R,\epsilon}$ form a neighbourhood basis for $\xi$, we can assume without loss of generality that $\eta = \xi$. We need to find a hyperplane $k$, such that $\xi \in U_{\go,k} \subset U_{R, \epsilon}(\xi)$.

Denote the geodesic representative of $\xi$ based at $\go$ by $b$. Since $\xi$ is $\kappa$-Morse, we find a family of excursion hyperplanes $(h_i)_i$ for $b$. Let $\eta \in \partial^{\kappa} X$ and $c$ the geodesic representative of $\eta$ based at $\go$. For any $n$, we denote $p_n := b \cap h_n$ and (if it exists) $q_n := c \cap h_n$. Suppose that $c$ crosses $h_{n+1}$. 
By Lemma \ref{lem:BoundofHyperplanesandGeodesics} there exists a constant $\dd$, which depends only on the excursion constant and $\kappa$, such that $$d(p_n,q_n) \leq \dd \kappa(t_n).$$

We have shown that for the geodesic representative $b$ of $\xi$ based at $\go$ and any geodesic ray $c$ based at $\go$ that intersects $h_{n+1}$, we can bound the distance $d(p_n,q_n)$. Suppose that $d(b(R), c(R)) \geq \epsilon$. Since distance functions are convex in a $\CAT$ space, this implies that, for all $t \geq R$, $d(b(t), c(t)) \geq \frac{t}{R} \epsilon$. Since $\kappa$ is sublinear, we find some $T$, such that for all $t \geq T$, $2\dd \kappa(t) < \frac{t}{R} \epsilon$. Define $N$ to be the first number such that $d(\go, p_N) \geq T$. We conclude that, for all $\eta \in U_{\go, h_{N+1}}$ with geodesic representative $c$,
\[ d(b(t_N), c(t_N)) \leq 2 d(p_N, q_N) \leq 2\dd \kappa(t_{N}) < \frac{t_{N}}{R} \epsilon. \]
The first inequality follows from the fact that $c$ can cross $h_N$ only at a time $s_N$ that satisfies $\vert t_N - s_N \vert \leq d(p_N, q_N)$ (because it is a geodesic and starts at the same point as the geodesic $b$). By the discussion on convexity above, this implies that $d(b(R), c(R)) < \epsilon$. Therefore, $U_{\go, h_{N+1}} \subset U_{R, \epsilon}(\xi)$. The lemma follows.
\end{proof}




\section{Continuity of the Roller map} \label{sec:RollerContinuity}

Let $X$ be a finite dimensional $\CAT$ cube complex and let $m$ be an upper bound for the degree of every vertex. Fix a base point $\go \in X^{(0)}$. Recall that we defined
\[ \partial_{R}^{\kappa} X = \{x \in \partial_R X| x \text{ can be represented by a } \kappa-\text{Morse combinatorial geodesic ray}\}. \]
Let $x \in \partial_{R}^{\kappa} X$ and let $\alpha$, $\beta$ be combinatorial geodesic rays starting at $\go$ that represent $x$. Suppose $\alpha$ is $\kappa$-Morse. Since $\alpha$ and $\beta$ both represent $x$ and start at the same vertex, we have $[\alpha \vert \beta]_{\go} = \infty$. By Lemma \ref{lem:GromovProductFellowTraveling}, this implies that $\alpha$ and $\beta$ are $\kappa$-fellow traveling and $\beta$ is $\kappa$-Morse as well. In particular, if we see $\alpha$ and $\beta$ as $\kappa$-Morse quasi-geodesics, they represent the same point in $\partial_{\kappa}X$ and thus a unique point in $\partial^\kappa_ \infty X$ using Lemma \ref{lem:injection}. This way, we obtain a map
\[ \Phi: \partial_R^{\kappa} X \rightarrow \partial^{\kappa}_\infty X. \]

We show that the map $\Phi$ is a continuous surjection.

\begin{theorem} \label{thm:Rollerprojection}
The natural map $\Phi: \partial_R^{\kappa} X \rightarrow \partial^{\kappa}_\infty X $ is a continuous surjection. For all $x, y \in \partial_R^{\kappa} X$, we have that $\Phi(x) = \Phi(y)$ if and only if $[x \vert y]_{\go} = \infty$. Furthermore, the induced quotient map $\overline{\Phi}:\overline{\partial_{R} ^\kappa X} \rightarrow \partial^{\kappa}_\infty X$ is a homeomorphism.
\end{theorem}

\begin{remark}
We can express the fibres of the map $\Phi$ using the notion of components in $\partial_R X$ and the partial order $\prec$ between components (see Definition \ref{def:RollerComponents} and the paragraphs after). The fibre $\Phi^{-1}(\xi)$ of any point $\xi \in \partial^{\kappa}_\infty X$ is the union of a component $C$ that is maximal with respect to $\prec$ and all components satisfying $C' \prec C$. In particular, we obtain that, whenever a component $C$ in $\partial_R X$ is represented by $\kappa$--Morse combinatorial geodesic rays, then $C$ can be contained in only one maximal element with respect to $\prec$.
\end{remark}

\begin{proof}

We first show surjectivity. Let $\xi \in \partial^{\kappa}_\infty X$ and $b$ the geodesic representative of $\xi$ based at $\go$. By Lemma \ref{lem:ApproximationofGeodesicRays}, there exists a combinatorial geodesic ray $\alpha$ based at $\go$ that crosses the same hyperplanes as $b$ and has finite Hausdorff distance to $b$. Let $x \in \partial_R^{\kappa} X$ be the point represented by $\alpha$. Clearly, $\Phi(x) = \xi$. Therefore, $\Phi$ is surjective.\\

Next, we show that $\Phi(x) = \Phi(y)$ if and only if $[x \vert y]_{\go} = \infty$. Let $x, y \in \partial_R^{\kappa} X$, suppose $[x \vert y]_{\go} = \infty$, and choose combinatorial representatives $\alpha, \beta$ that start at $\go$. Since $\alpha$ has to be $\kappa$-Morse, Lemma \ref{lem:GromovProductFellowTraveling} implies that $\alpha$ and $\beta$ $\kappa$-fellow travel and, therefore, $\Phi(x) = \Phi(y)$.

Now suppose, $[x \vert y]_{\go} < \infty$. Choose a combinatorial geodesic from $\go$ to the median $m(\go, x, y)$ and combinatorial representatives of $x$ and $y$ respectively that start at $m(\go,x,y)$. Concatenating the geodesic segment with one of the combinatorial geodesic rays yields representatives of $x$ and $y$ respectively that start at $\go$ and are identical until they have crossed all hyperplanes in $\mathcal{W}(\go \vert x,y)$. Denote these two representatives by $\alpha, \beta$. We claim that, after separating, $\alpha$ and $\beta$ diverge linearly. Let $T$ be the time when $\alpha(T) = m(\go, x, y) = \beta(T)$. Since every hyperplane crossed by $\alpha$ after time $T$ separates the unbounded part of $\alpha$ from $\beta$, we have
\[ d(\alpha(T + n), \beta) \geq \frac{n}{\sqrt{m}}, \]
which grows linearly in $n$. We conclude that $\alpha, \beta$ are not $\kappa$-fellow traveling. Therefore, $\Phi(x) \neq \Phi(y)$, which implies that the fibres of $\Phi$ are as described above.\\

We now prove continuity of $\Phi$. Let $\xi \in \partial^{\kappa}_\infty X$ and $b$ the geodesic representative of $\xi$ based at $\go$. By Lemma \ref{lem:ApproximationofGeodesicRays}, there exists a combinatorial geodesic ray $\alpha$ based at $\go$ that crosses the same hyperplanes as $b$ and which has finite Hausdorff distance to $b$. Denote the point on the Roller boundary represented by $\alpha$ as $x$. Clearly, $x \in \Phi^{-1}(\xi)$.

By Theorem \ref{thm: characterizing sublinear Morse geodesics using hyperplanes} and Theorem \ref{thm:Morse iff contracting iff excursion}, there exists a family of excursion hyperplanes $(h_i)_i$ for $b$. Consider the sets
\[ U_{\go, h_i} = \{ \eta \in \partial_\infty^{\kappa} X \vert \eta_{\go} \text{ crosses } h_i \}, \]
where $\eta_{\go}$ denotes the unique geodesic representative of $\eta$ that starts at $\go$. By Lemma \ref{lem:BasisforHYP}, including its proof, and Theorem \ref{thm:VisualequalsHyp} the family $\{ U_{\go, h_i} \}_i$ forms a neighbourhood basis of the cone topology at $\xi$. In order to prove continuity of $\Phi$, we are left to prove that for every $i$, we find a hyperplane $k$ such that the set
\[ V_{\go,k} := \{ y \in \partial^\kappa_ R X \vert k \in \mathcal{W}(\go \vert y) \}, \]
which is open in $\partial_R^{\kappa} X$ with the Roller topology, satisfies
\[ \Phi^{-1}(\xi) \subset V_{\go, k}, \]
\[ \Phi(V_{\go,k}) \subset U_{\go,h_i}. \]

We claim that, for all $i$, $\Phi(V_{\go,h_{i+1}}) \subset U_{\go, h_i}$ and $\Phi^{-1}(\xi) \subset V_{\go, h_{i+1}}$. We first show that for all $x' \in \Phi^{-1}(\xi)$, any combinatorial representative $\alpha$ of $x'$ starting at $\go$ has to cross all excursion hyperplanes of $b$. If $\alpha$ does not cross an excursion hyperplane $h$ of $b$, then we have
\[ d(b(t), \alpha) \geq d(b(t), h^{-}), \]
where $h^{-}$ denotes the halfspace of $h$ containing $\go$ and thus $\alpha$. Since $b$ crosses $h$ the function $d(b(t), h^{-})$ is at least linearly increasing (eventually), implying that $d(b(t), \alpha)$ grows at least linearly. This contradicts the fact that $\Phi(x') = \xi$ and therefore, $\Phi^{-1}(\xi) \subset V_{\go, h_{i+1}}$.

We are left to prove that $\Phi(V_{\go, h_{i+1}}) \subset U_{\go, h_i}$. Let $y \in V_{\go, h_{i+1}}$ and let $\tilde{\beta}$ be a combinatorial geodesic ray, representing $y$, starting at $\go$. Since $y \in \partial^\kappa_R X$, $\tilde{\beta}$ is $\kappa$-Morse. Therefore, there exists a $\kappa$-Morse geodesic ray $c$, starting at $\go$, that is $\kappa$-fellow traveling with $\tilde{\beta}$ and a combinatorial geodesic ray $\beta$ that crosses exactly the same hyperplanes as $c$. Suppose that $c$, and therefore $\beta$, does not cross $h_i$. We will show that this implies that $\Phi([\beta]) \neq \Phi([\tilde{\beta}])$, which is a contradiction, as both of them are $\kappa$-fellow traveling along $c$.

Since $\beta$ does not cross $h_i$ it cannot cross $h_{i+1}$ either. Therefore, every $k \in \mathcal{W}(\beta) \cap \mathcal{W}(\tilde{\beta})$ either crosses $h_{i+1}$ or it crosses $\tilde{\beta}$ before $\tilde{\beta}$ crosses $h_{i+1}$. Since $h_i$ and $h_{i+1}$ are well-separated and the segment of $\tilde{\beta}$ before crossing $h_{i+1}$ is finite, we conclude that $\mathcal{W}(\beta) \cap \mathcal{W}(\tilde{\beta})$ is a finite set. In other words, their combinatorial Gromov product satisfies $\left[ [ \beta ] \vert [ \tilde{\beta} ] \right]_{\go} < \infty$. From the characterisation of the fibres of $\Phi$, we conclude that $\Phi([\beta]) \neq \Phi([\tilde{\beta}])$, a contradiction. This implies that $\beta$, and thus $c$ crosses $h_i$. Therefore, $\Phi(y) = \Phi([\tilde{\beta}]) \in U_{\go, h_i}$, which implies continuity of $\Phi$.\\

Saying that two points in $\partial_{R}^{\kappa} X$ are equivalent if and only if they lie in the same fibre of $\Phi$, we obtain a quotient space, denoted $\overline{\partial_{R}^{\kappa} X}$, and a bijective, continuous map $\overline{\Phi} : \overline{ \partial_{R}^{\kappa} X} \rightarrow \partial^{\kappa}_\infty X$. To prove the Theorem, we are left to show that $\overline{\Phi}^{-1}$ is continuous. We begin by proving the following result: Given points $x_n, x \in \partial_{R}^{\kappa} X$ and $y \in \partial_R X$ such that $x_n \rightarrow y$ and $\Phi(x_n) \rightarrow \Phi(x)$, we have $[x \vert y]_{\go} = \infty$.

Let $b_n, b$ be geodesic representatives of $\Phi(x_n)$ and $\Phi(x)$ respectively, all starting at $\go$. By Lemma \ref{lem:ApproximationofGeodesicRays}, there exist combinatorial geodesics $\alpha_n$ and $\alpha$ based at $\go$ that are asymptotic to and cross the same hyperplanes as $b_n$ and $b$ respectively. Denote the points in $\partial_R X$ represented by $\alpha_n$ and $\alpha$ by $z_n$ and $z$ respectively. By construction, $\Phi(x_n) = \Phi(z_n)$ and $\Phi(x) = \Phi(z)$ and, therefore, $[x_n \vert z_n]_{\go} = [x \vert z]_{\go} = \infty$. By Theorem \ref{thm: characterizing sublinear Morse geodesics using hyperplanes} and Theorem \ref{thm:Morse iff contracting iff excursion}, we find a sequence of excursion hyperplanes $(h_i)_i$ for $b$, which are thus crossed by $\alpha$. Let $\alpha_x$ and $\alpha_y$ be combinatorial geodesic rays based at $\go$ that represent $x$ and $y$ respectively. We claim that both $\alpha_x$ and $\alpha_y$ cross $h_i$ for all $i$.

Suppose $\alpha_x$ did not cross $h_i$. Since $[x \vert z]_{\go} = \infty$, there are infinitely many hyperplanes crossed by both $\alpha$ and $\alpha_x$ after $\alpha$ crosses $h_{i+1}$. All these hyperplanes have to intersect both $h_i$ and $h_{i+1}$, which contradicts the fact that $h_i$ and $h_{i+1}$ are $\kappa$-well-separated. Therefore, $\alpha_x$ crosses $h_i$.

Suppose $\alpha_y$ did not cross $h_i$. Let $\alpha_{x_n}$ be a combinatorial geodesic ray, starting at $\go$ representing $x_n$. Since we assume that $x_n \rightarrow y$ in Roller topology and $\Phi(x_n) \rightarrow \Phi(x)$ in cone topology (which equals $\Hyp$), we find a number $N$ such that for all $n \geq N, \alpha_{x_n}$ does not cross $h_i$ and $b_n$ crosses $h_{i+1}$. By the same argument as above, there can be at most finitely many hyperplanes crossed by both $\alpha_{x_n}$ and $b_n$ after $b_n$ crosses $h_{i+1}$. Since $b_n$ crosses $h_{i+1}$ after finite time, the total number of hyperplanes crossed by both $\alpha_{x_n}$ and $b_n$ has to be finite. Since $b_n$ and $\alpha_n$ cross the same hyperplanes, we conclude that $[x_n \vert z_n]_{\go} < \infty$, which is a contradiction. Therefore, $\alpha_y$ crosses $h_i$.

We conclude that for all $i$, $h_i$ separates $\go$ from both $x$ and $y$. Therefore, $[x \vert y ]_{\go} = \infty$.

We now prove that $\overline{\Phi}^{-1}$ is sequentially continuous. We denote the projection map onto the quotient by $p : \partial^\kappa_R X \rightarrow \overline{ \partial^\kappa_R X}$. Let $x_n, x \in \partial^\kappa_R X$ such that $\Phi(x_n) \rightarrow \Phi(x)$. Suppose that $p(x_n) \nrightarrow p(x)$. Then there exists an open neighbourhood $V \subset \overline{\partial^\kappa_R X}$ of $p(x)$ and a subsequence of $x_n$, again denoted by $x_n$, such that $x_n \notin p^{-1}(V)$ for all $n$. Since $X$ is locally finite, $\partial_R X$ is compact and, by moving to a subsequence if necessary, we can assume that $x_n \rightarrow y \in \partial_R X$. We have shown above that $[ x \vert y]_{\go} = \infty$. By Lemma \ref{lem:GromovProductFellowTraveling}, we see that $y \in \partial^\kappa_R X$ and $[x \vert y]_{\go} = \infty$. By our characterisation of the fibres of $\Phi$, we conclude that $\Phi(y) = \Phi(x)$ and $p(y) = p(x) \in V$. Therefore, $\lim_{n \rightarrow \infty} x_n = y \in p^{-1}(V)$, contradicting the assumption on $x_n$. We conclude that $p(x_n) \rightarrow p(x)$. Therefore, $\overline{\Phi}^{-1}$ is sequentially continuous and, since $\partial^{\kappa}_\infty X$ is first countable, it is continuous. This completes the proof.
\end{proof}




\section{Embedding the sublinear boundary into a Gromov boundary} \label{sec:EmbeddingintoGromovboundary}

The main goal of this section is to prove Theorem \ref{introthm:cont injection}. 

\begin{definition}[Parallel subcomplexes]\label{parallel} Let $X$ be a finite-dimensional $\CAT$ cube complex. Two convex subcomplexes $H_1, H_2$ are said to be \emph{parallel} if for any hyperplane $h,$ we have $h \cap H_1 \neq \emptyset$ if and only if $h \cap H_2 \neq \emptyset$.
\end{definition}

\begin{definition}[Factor systems \cite{Behrstock2017}] \label{def:factorsystem}Let $X$ be a finite-dimensional $\CAT$ cube complex. A \emph{factor system}, denoted $\mathfrak{F}$, is a collection of subcomplexes 
of $X$ such that:

\begin{enumerate}
    \item $X \in \mathfrak{F}$.
    \item Each $F \in \mathfrak{F}$ is a nonempty convex subcomplex of $X$
    \item There exists $\cc_{1} \geq 1$ such that for all $x \in X^{(0)}$ at most $\cc_{1}$ elements of $\mathfrak{F}$ contain $x$.
    \item Every nontrivial convex subcomplex parallel to a hyperplane of $X$ is in $\mathfrak{F}$.
    \item There exists $\cc_{2}$ such that for all $F, F' \in \mathfrak{F}$, either $\sg_{F}(F') \in \mathfrak{F}$ or $\diam(\sg_{F}(F')) \leq \cc_{2}$.
\end{enumerate}

\end{definition}

\begin{remark}
All known cocompact CAT(0) cube complexes admit factor systems \cite{HagenSusse}.
\end{remark}

We also use the following generalization of the contact graph introduced by Genevois in \cite{HypInCube}: 

\begin{definition}[The well-separation space] \label{def:Morse-detecting space}
Let $X$ be a finite-dimensional $\CAT$ cube complex. For each non-negative integer $k$, the \emph{$k$-well-separation space}, denoted by $Y_k$, is defined to be the set whose elements are the vertices of $X$, with the following distance function. For $x,y \in X^{(0)}$, the $k$-distance between $x,y$, denoted by $d_k(x,y)$, is defined to be the cardinality of a maximal collection of $k$-well-separated hyperplanes separating $x,y$. We will refer to the space $(Y_k, d_k)$ as the $k$-well-separation space.

\end{definition}

It's worth remarking that when $k=0$, the $0$-well-separation space is quasi-isometric to the contact graph (Fact 6.50 of \cite{HypInCube}). Furthermore, we have the following:

\begin{proposition}[{\cite[Proposition 6.54]{HypInCube}}]\label{prop: Gen's metric is hyperbolic} Let $X$ be a finite-dimensional $\CAT$ cube complex and let $Y_k$ be the $k$-well-separation space as in Definition \ref{def:Morse-detecting space}. For any non-negative integer $k$, the metric space $(Y_k,d_k)$ is $9(k+2)$ hyperbolic.
\end{proposition}

\begin{remark} \label{rmk: projection to vertices}

For each $x \in X,$ let $n$ be the largest integer so that $x \in [0,1]^n$. We define $v(x)$ to be a $0$-cell (a vertex) of the cube $[0,1]^n$. For any $x,y \in X$, we denote the cardinality of the maximal collection of $k$-well-separated hyperplanes separating $x,y$ by $l_k(x,y)$. Since $X$ is a finite dimensional cube complex, there exists a constant $\cc$ depending only on the dimension of $X$ such that the following holds:

\begin{itemize}

    \item We have: $ d_k(v(x),v(y))-\cc\leq l_k(x,y) \leq d_k(v(x),v(y))+\cc. $

\item We have $ d(x,y)-\cc\leq d^{(1)}(v(x),v(y)) \leq d(x,y)+\cc,$ where $d^{(1)}$ and $d$ are the combinatorial and $\CAT$ distance respectively.

\item We have $d_k(v(x), v(y)) \leq d(x,y)+\cc,$ where $d_k$ and $d$ are the distances in $Y_k$ and $X$ respectively.

\end{itemize}

\end{remark}

\begin{lemma} \label{lem:ProjectingtoQuasiGeodesics}
Let $X$ be a $\CAT$ cube complex and let $(Y_k,d_k)$ be the $k$-well-separation space. Every ($\CAT$ or combinatorial) geodesic in $X$ projects to a uniform unparameterized quasi-geodesic in $(Y_k,d_k)$. More precisely, we have the following:

\begin{enumerate}
    \item For a combinatorial geodesic $\alpha$ in $X^{(1)}$, if $x,y,z \in X^{(0)}$ are points with $x=\alpha(t_1), y=\alpha(t_2)$ and $z=\alpha(t_3)$ with $t_1 \leq t_2 \leq t_3,$ we have $d_k(x, z) \geq d_k(x, y)+d_k(y, z)-k-3$. 
    
    \item For a $\CAT$ geodesic $c$ in $X$, if $x,y,z \in X$ are points with $x=c(t_1), y=c(t_2)$ and $z=c(t_3)$ with $t_1 \leq t_2 \leq t_3,$ we have $d_k(x, z) \geq d_k(x, y)+d_k(y, z)-k-\cc-3$, where $\cc$ is a constant depending only on the dimension of $X$ (and not on the geodesic $c$).
\end{enumerate}

\end{lemma}

\begin{proof} We only prove the first part of the assertion, the second part follows immediately by Remark \ref{rmk: projection to vertices}.
Let $L_1, L_2$ denote $d_k(x, y), d_k(y, z)$ respectively. By definition, $L_1, L_2$ are the maximal cardinality of $k$-well-separated hyperplanes separating $x,y$ and $y,z$ respectively. In other words, if $\mathcal{W}_1, \mathcal{W}_2$ denote maximal sets of $k$-well-separated hyperplanes separating $x, y$ and $y,z$ respectively, then $|\mathcal{W}_1|=L_1$ and $|\mathcal{W}_2|=L_2$. Let $\mathcal{W}_1=\{h_1,h_2, \dots h_{L_1} \}$ and define $t_i$ such that  $\alpha(t_i)=h_i \cap \alpha$ and $t_i <t_{i+1}$ for all $1 \leq i \leq L_1 $. Similarly, let  $\mathcal{W}_2=\{h'_1,h'_2, \dots h'_{L_2}\}$ where  $\alpha(s_i)=h_i \cap \alpha$ and $s_i <s_{i+1}$ for all $1 \leq i \leq L_2 $.

\begin{itemize}
    \item Case 1: If $h_{L_1}$ and $h'_2$ are disjoint, then any two hyperplanes $h_i \in \mathcal{C}_1$ and $h'_{j} \in \mathcal{C}_2 -\{h'_1\}$ must be $k$-well-separated since every hyperplane $h$ meeting both $h_i,h'_{j}$ must also meet $h'_{1}$. Therefore, every such $h$ must meet the hyperplanes $h'_1, h'_j$, but since $h'_1, h'_j$ are $k$-well-separated, there can be at most $k$ such hyperplanes containing no facing triple. That is to say, $d_k(x,z) \geq d_k(x, y)+d_k(y, z)-1$.
    
    \item Case 2: (See figure \ref{fig:ProjectiontoUnparametrisedGeodesics}.) If $h_{L_1}$ and $h'_2$ are not disjoint, choose the smallest $m$ so that $h_{m} \cap h'_2 \neq \emptyset.$ Notice that, since $h'_1,h'_2$ are disjoint, we have $h_m \cap h'_1 \neq \emptyset $. In particular, the hyperplanes $h'_1,h'_2$ must both cross $h_m,h_{m+1},\dots, h_{L_1}.$ Since $h_m,h_{m+1},\dots, h_{L_1}$ all meet the geodesic $\alpha$, they can't contain a facing triple and as $h'_1,h'_2$ are $k$-well-separated, we have $|\{h_m,h_{m+1},\dots, h_{L_1}\}| \leq k.$ That is to say, $m \geq L_1-k$. In particular, $h'_2$ must be disjoint from the hyperplanes $h_1,h_2,\dots, h_{L_1-k-1}.$ Hence, proceeding as in argument for Case 1, the hyperplanes $h_1,h_2,\dots, h_{L_1-k-1},h'_3,h'_4, \dots h'_{L_2}$ are pairwise $k$-well-separated. In other words, $d_k(x,z) \geq (L_1-k-1) +(L_2-2)=L_1+L_2-k-3$. Notice that the projection is a uniform unparameterized quasi-geodesic as the constant  $k+3$ is independent of the projected geodesic $\alpha.$
\end{itemize}


\end{proof}

We conclude from Lemma \ref{lem:ProjectingtoQuasiGeodesics} that, if a geodesic ray crosses an infinite sequence of $k$-well-separated hyperplanes, then it projects to an infinite, unparametrised quasi-geodesic and thus defines a point in the Gromov boundary of $(Y_k,d_k)$. More precisely, we have the following corollary.

\begin{corollary}\label{corollary: crossing infinite sequence of well-separated hyperplanes gives a point in the Gromov boundary}
Every ($\CAT$ or combinatorial) geodesic ray crossing an infinite sequence of $k$-well-separated hyperplanes projects to an infinite unparameterized quasi-geodesic in $(Y_k,d_k)$ with uniform quasi-geodesic constants (independent of the geodesic) defining a point in  $\partial Y_k$, the Gromov boundary of $Y_k.$
\end{corollary}

The following lemma appears in the first version of \cite{HypInCube} on the arXiv, however, due to technical issues, it was omitted from the second version. With Genevois's blessing, the statement and the argument were reproduced in \cite{Murray-Qing-Zalloum}.

\begin{lemma} [{\cite[Lemma 4.26]{Murray-Qing-Zalloum}}]\label{lem: well-separated implies L-well-separated}
Let $X$ be a cocompact $\CAT$ cube complex with a factor system. There exists a constant $L \geq 0$ such that any two hyperplanes of $X$
 are either $L$-well-separated or are not well-separated.
\end{lemma}

The constant from the above lemma will be referred to frequently in this paper, therefore, we give it a name:

\begin{definition}(the separation constant) Let $X$ be cocompact $\CAT$ cube complex with a factor system. The  \emph{the separation constant} of $X$ is the smallest integer $L \geq 0$ such that any two hyperplanes of $X$ are either $L$-well-separated or are not well-separated. Such a constant exists using Lemma \ref{lem: well-separated implies L-well-separated}.

\end{definition}

 In light of Lemma \ref{lem: well-separated implies L-well-separated}, Corollary \ref{corollary:  crossing infinite sequence of well-separated hyperplanes gives a point in the Gromov boundary}, Theorem \ref{thm: characterizing sublinear Morse geodesics using hyperplanes}, and Theorem \ref{thm:Morse iff contracting iff excursion} we get the following.

\begin{corollary} 
Let $X$ be a cocompact $\CAT$ cube complex with a factor system and let $L$ be the separation constant. Every infinite $\kappa$-Morse ($\CAT$ or combinatorial) geodesic projects to an infinite unparameterized quasi-geodesic in $(Y_L,d_L)$ with uniform quasi-geodesic constants (independent of the projected geodesic) defining a point in $\partial Y_L$, the Gromov boundary of $Y_L.$
\begin{proof}
The proof follows immediately from Lemma \ref{lem: well-separated implies L-well-separated}, Corollary \ref{corollary:  crossing infinite sequence of well-separated hyperplanes gives a point in the Gromov boundary}, Theorem \ref{thm: characterizing sublinear Morse geodesics using hyperplanes} and Theorem \ref{thm:Morse iff contracting iff excursion}.
\end{proof}

\end{corollary}

The previous corollary implies the existence of a map from the $\kappa$-Morse boundary of a $\CAT$ cube complex to the Gromov boundary of the hyperbolic space $Y_L$. The following lemma is to show that this map is injective.

\begin{lemma}[Injectivity of the map] \label{lem: injectivity of the map}
Let $X$ be a cocompact $\CAT$ cube complex with a factor system and let $L$ be the separation constant. Let $c,c'$ be two distinct $\CAT$ $\kappa$-Morse geodesic rays with $c(0)=c'(0)=\go$. For each $m$, there exists $t$ large enough so that the number of $L$-well-separated hyperplanes separating $c(t)$ from $c'$ is at least $m.$
\end{lemma}

\begin{proof} Notice that since $c,c'$ are distinct, neither of them can be in the $\kappa$-neighborhood of the other. Let $\mathcal{H}=\{h_i\}$ denote the excursion hyperplane's for $c$  and let $c(t_i)$ denote $h_i \cap c$. Suppose that the set $\mathcal{H}$ is ordered based on the order in which the excursion hyperplanes cross $c$. Using Lemma \ref{lem:BoundofHyperplanesandGeodesics}, only finitely many hyperplanes of $\mathcal{H}$ can meet $c'$. Let $j$ be the largest integer so that $h_{j}$ crosses  $c'$. Define $K=d(\go, h_j \cap c)$. Notice that for $t>K,$ if $\mathcal{W}_t$ denotes the collection of $L$-well-separated excursion hyperplanes separating $c(t_j)$ from $c(t)$, then $|\mathcal{W}_t|\rightarrow \infty$ as $t \rightarrow \infty.$ For a fixed $t>K$, since hyperplanes of $\mathcal{W}_t$ can't cross $c'$, they must cross every geodesic starting at $c(t)$ and ending on $c'$. This implies that for each $m$, there exists a $t$ such that the number of $L$-well-separated hyperplanes spearating $c(t)$ from $c'$ is at least $m$. This gives the desired conclusion.
\end{proof}

The well-separation metric space given in Definition \ref{def:Morse-detecting space} is not a geodesic metric space. We will use the well-separation space defined above to define a geodesic metric space with a comparable distance function. More precisely, for an integer $k$, let $\Gamma_k$ be the graph whose vertices are $X^{(0)}$ with two vertices $x,y \in X^{(0)}$ connected by an edge if and only if the number of $k$-well-separated hyperplanes separating $x,y$ is at most $10k+4$.

\begin{definition}[The well-separation graph] \label{def:Morse-detecting graph}
Let $X$ be a finite dimensional $\CAT$ cube complex. For each integer $k$, the \emph{$k$-well-separation graph}, denoted by $\Gamma_k$ is defined to be the graph whose vertices are the vertices of $X$ with an edge connecting two distinct vertices $x,y$ if the number of $k$-well-separated hyperplanes separating $x$ from $y$ is bounded from above by $10k+4.$

\end{definition}

We now show that the spaces $(Y_k,d_k)$ and $(\Gamma_k,d_{\Gamma_k})$ are bilipschitz equivalent to each other. More precisely, we show the following.

\begin{lemma}[Bilipschitz equivalence] \label{lem:updating to a graph}
Let $X$ be a finite-dimensional $\CAT$ cube complex and let $(Y_k,d_k), (\Gamma_k,d_{\Gamma_k})$ be the $k$-well-separations space and graph respectively. For any two points $x,y \in X^{(0)},$ we have $\frac{1}{10k+4}d_k(x,y) \leq d_{\Gamma_k}(x,y) \leq d_k(x,y)$.
\end{lemma}

\begin{proof}
We first show that $\frac{1}{10k+4}d_k \leq d_{\Gamma_k}$. For $x,y \in X^{(0)}$, let $p$ be any edge path in $\Gamma$ connecting $x,y$. Denote the vertices and the edges of this path by $x_0=x,x_1,x_2, \dots x_m=y$ and $e_1,e_2, \dots, e_m$, respectively. Let $\Tilde{p}$ denote some edge path in $X^{(1)}$ connecting $x_0=x,x_1,x_2, \dots x_m=y$ such that the subpath $\tilde{p}_i$ of $\tilde{p}$ connecting $x_{i-1}$ to $x_i$ is a combinatorial geodesic for all $i$. Let $\mathcal{W}=\{h_1,h_2, \dots,h_n\}$ denote a maximal collection of $k$-well-separated hyperplanes separating $x,y.$ Every such hyperplane must meet $\Tilde{p}$, and since $x_i,x_{i+1}$ are connected by an edge $e_i$ in $\Gamma_k,$ each subpath $\tilde{p}_i$ can be crossed by at most $10k+4$ hyperplanes of $\mathcal{C}.$ Therefore, $\frac{n}{10k+4} \leq m$ which gives that $\frac{1}{10k+4}d_k(x,y) \leq d_{\Gamma_k}(x,y)$.

Now we show that $d_{\Gamma_k}(x,y) \leq d_k(x,y)$. In order to do so, we will construct an edge path in $\Gamma$ of length $d_k(x,y).$ Let $\alpha$ be a combinatorial geodesic in $X$ connecting $x,y$ and let $\mathcal{W}=\{h_1,h_2, \dots,h_n\}$ be a maximal collection of $k$-well-separated hyperplanes separating $x,y$. Let $x_0=x,$ $x_i=\alpha(t_i)$ denote the vertex on $\alpha$ immediately after crossing the hyperplane $h_i$ for all $1 \leq i \leq n-1$. We claim that there is an edge in $\Gamma$ connecting $x_i$ to $x_{i+1}$. In other words, we claim that the number of $k$-well-separated hyperplanes separating $x_i,x_{i+1}$ is at most $10k+4.$ Suppose not, that is, suppose that $x_i, x_{i+1}$ are separated by a collection of $k$-well-separated hyperplanes $\mathcal{W}'=\{h_1',h_2', \dots, h'_{m}\}$ with $m \geq 10k+5$. Notice that $h_{i+1} \notin \mathcal{W}'$ as if $h_{i+1} \in  \mathcal{W}'$, then, since the hyperplanes $h_{i-1},h_i$ are $k$-well-separated, at most $k$ hyperplanes of the collection $\mathcal{W}'$ can cross them. Therefore, the collection $(\mathcal{W}-\{h_i,h_{i-1}\}) \cup \{h_i' \in \mathcal{W}'| h_i' \text{  does not cross both }h_i,h_{i-1}\}$ is a collection of $k$-well-separated hyperplanes separating $x,y,$ whose cardinality is at least $(n-2)+(10k+5-k)>n$  contradicting the maximality assumption. Now, notice that for any five hyperplanes $h_i',h_j', h_l', h'_o, h'_p \in \mathcal{W'}$ ordered based on the order they are crossed by $\alpha,$ at least one of them must cross  either both $h_i,h_{i-1}$ or both $h_{i+1},h_{i+2}$ as otherwise one gets a larger collection of $k$-well-separated hyperplanes separating $x,y.$ In other words, if none of the five hyperplanes  $h_i',h_j', h_l', h'_o, h'_p$ cross both $h_i,h_{i-1}$ or both $h_{i+1},h_{i+2}$, then, the set consisting of $(\mathcal{W}-\{h_i, h_{i+1}\}) \cup \{h_j',h_l',h_o'\} $ is a set of $k$-well-separated hyperplanes separating $x,y$ with cardinality $n+1$ contradicting the maximality assumption, see Figure \ref{fig:BilipschitzEquivalence}. Notice that $10k+5$ is exactly $5(2k+1)$, in other words, grouping elements of $\mathcal{W'}$ into subsets of cardinality exactly 5 yields exactly $2k+1$ such subsets. Each of these subsets contains at least one hyperplane crossing either both $h_i,h_{i-1}$ or both $h_{i+1},h_{i+2}$. This implies that either $h_i,h_{i-1}$ or $h_{i+1},h_{i+2}$ is crossed by $k+1$ hyperplanes, which is a contradiction. The above argument shows that there is an edge connecting $x_i$ to $x_{i+1}$ for all $ 2\leq i \leq m-2$. It remains to show that there are three edges $e_1, e_2,e_3$ connecting $\{x,x_1\}$, $\{x_1,x_2\},\{x_{m-1},y\}$ respectively, but an almost identical argument to the above shows that.

\begin{figure}
\begin{tikzpicture}[scale=1.50]
\draw [] (-3,1) -- (3,1);
\node [right] at (3,1) {$\alpha$};
\draw (-2.8, -0.5) -- (-2.8, 1.5);
\draw (-1.4, -0.5) -- (-1.4, 1.5);
\draw (0.8, -0.5) -- (0.8, 1.5);
\draw (2, -0.5) -- (2, 1.5);
\node [above] at (-2.8, 1.5) {$h_{i-1}$};
\node [above] at (-1.4, 1.5) {$h_{i}$};
\node [above] at (0.8, 1.5) {$h_{i+1}$};
\node [above] at (2, 1.5) {$h_{i+2}$};
\draw [fill] (-1.2, 1) circle [radius=0.02cm];
\node [below] at (-1.2, 1) {$x_i$};
\draw [fill] (1, 1) circle [radius=0.02cm];
\node [below] at (1.1, 1) {$x_{i+1}$};
\draw (-0.5, 1.5) to [out = -90, in = 0] (-1.4, 0.7);
\draw (-1.4, 0.7) to [out = 180, in = 90] (-2.4, -0.5);
\node [below] at (-2.4, -0.5) {$h'_i$};
\draw (-0.4, 1.5) to [out = -90, in = 0] (-1.4, 0.5);
\draw (-1.4, 0.5) to [out = 180, in = 90] (-2.1, -0.5);
\node [below] at (-2.1, -0.5) {$h'_j$};
\draw (-0.3, 1.5) to [out = -90, in = 0] (-1.4, 0.3);
\draw (-1.4, 0.3) to [out = 180, in = 90] (-1.8, -0.5);
\node [below] at (-1.8, -0.5) {$h'_l$};
\draw (-0.2, 1.5) to [out = -90, in = 180] (0.8, 0.4);
\draw (0.8, 0.4) to [out = 0, in = 90] (1.2, -0.5);
\node [below] at (1.2, -0.5) {$h'_o$};
\draw (-0.1, 1.5) to [out = -90, in = 180] (0.8, 0.6);
\draw (0.8, 0.6) to [out = 0, in = 90] (1.6, -0.5);
\node [below] at (1.6, -0.5) {$h'_p$};
\end{tikzpicture}
\caption{Given five hyperplanes in $\mathcal{W}'$ that intersect neither $h_{i-1}$ nor $h_{i+2}$, we can take the ``middle" three and obtain that $\mathcal{W} \setminus \{ h_i, h_{i+1} \} \cup \{ h'_j, h'_l, h'_o \}$ is a set of $k$-well-separated hyperplanes separating $x$ from $y$ that is strictly larger than $d_k(x, y)$.}
\label{fig:BilipschitzEquivalence}
\end{figure}
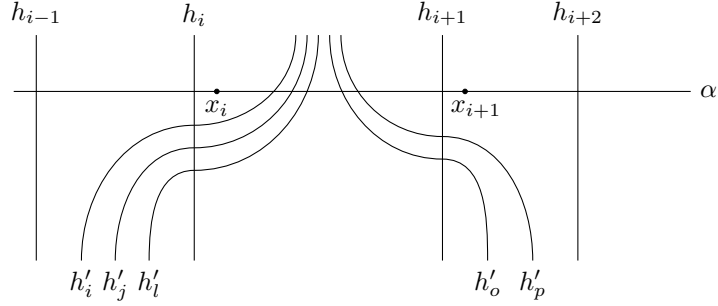

\end{proof}

\begin{corollary} \label{cor: new vertex map}
Let $X$ be finite dimensional $\CAT$ cube complex $X$ and let $k$ be a positive integer. There exists a constant $\cc$ and a map $v:(X,d) \rightarrow (\Gamma_k, d_{\Gamma_k})$ such that $d_{\Gamma_k}(v(x), v(y)) \leq d_k(v(x), v(y))  \leq d(x,y) +\cc.$
\end{corollary}

\begin{proof}
This follows using Lemma \ref{lem:updating to a graph} and part 3 of Remark \ref{rmk: projection to vertices}.
\end{proof}

\begin{corollary} \label{cor: sublinearly Morse geodesics define boundary points}
Let $X$ be a cocompact $\CAT$ cube complex with a factor system and let $L$ be the separation constant. Every infinite $\kappa$-Morse ($\CAT$ or combinatorial) geodesic projects to an infinite unparameterized quasi-geodesic in $(\Gamma_L,d_{\Gamma_L})$ with uniform quasi-geodesic constants (independent of the projected geodesic) defining a point in $\partial \Gamma_L$, the Gromov boundary of $\Gamma_L.$

\end{corollary}

\begin{proof}
The proof follows immediately from Lemma \ref{lem: well-separated implies L-well-separated}, Corollary \ref{corollary:  crossing infinite sequence of well-separated hyperplanes gives a point in the Gromov boundary}, Theorem \ref{thm: characterizing sublinear Morse geodesics using hyperplanes} and Theorem \ref{thm:Morse iff contracting iff excursion}.
\end{proof}

\begin{definition}[coarsely distance-decreasing] A map $f:(X,d) \rightarrow (Y,d')$ is said to be coarsely distance-decreasing if there exists constants $K,C \geq 0$ such that $d(f(x), f(y)) \leq K d(x,y)+C$ for all $x,y \in X.$

\end{definition}

Recall from Remark \ref{rmk: projection to vertices} that for a $\CAT$ cube complex $X$, there is a coarsely distance decreasing map $v: X \rightarrow X^{(0)}$ and a constant $\cc$ such that for any $x,y \in X$, if $l_k(x,y)$ denotes the cardinality of a maximal collection of $k$-well-separated hyperplanes separating $x,y$, then

\begin{itemize}

    \item We have: $ d_k(v(x),v(y))-\cc\leq l_k(x,y) \leq d_k(v(x),v(y))+\cc. $

\item We have $ d(x,y)-\cc\leq d^{(1)}(v(x),v(y)) \leq d(x,y)+\cc,$ where $d^{(1)}$ and $d$ are the combinatorial and $\CAT$ distance respectively.

\end{itemize}

The following lemma is phrased in slightly more general terms than needed, but it will immediately imply that $\partial_{\infty}^{\kappa} X$ continuously injects in the Gromov boundary of the geodesic hyperbolic metric space $\Gamma_L$. 

 Recall that convergence in the Gromov boundary of a $\delta$-hyperbolic space is characterized as follows. If $l \geq 2\delta+1,$ and $c_n,c$ are geodesic rays with $c_n(0)=c(0)$, then.:\\

$$[c_n] \rightarrow [c] \text{ iff } \forall t, \exists k=k_t \text{ such that }  d(c_n(t), c(t)) \leq l,\, \forall n \geq k, $$

\begin{lemma} \label{lem: continuous map gnrl settings}
Let $X$ be a \CAT \text { }space, $\Gamma$ a geodesic hyperbolic space and $p: X \rightarrow \Gamma$ a coarsely distance decreasing map. Suppose that $B \subseteq \partial _{\infty} X$ is a subspace of the visual boundary such that every element of $B$ is represented by a geodesic ray starting at $\go$ which projects to an infinite unparameterized quasi-geodesics with uniform constants in $\Gamma,$ then $p$ induces a continuous map $i_p: B \rightarrow \partial \Gamma.$
\end{lemma}

\begin{proof} Let $p: X \rightarrow \Gamma$ be as in the statement of the lemma. In particular, there exists constants $K, C$ such that $d(p(x),p(y)) \leq K d(x,y)+C.$
Suppose that $B$ is a subspace as in the statement of the lemma (in particular, it is equipped with the subspace topology of the visual boundary) and let $\cc_1,\cc_2$ be the constants so that $p(B)$ are all $(\cc_1, \cc_2)$-quasi-geodesics, after reparametrizing. Fix $\go \in X$ such that every geodesic ray in $B$ starts at $\go.$ Let $b_n \rightarrow b \in B$, and let $q_n, q$ denote the infinite $(\cc_1, \cc_2)$-quasi-geodesics which are the images of $b_n, b$ respectively. Let $c_n, c$ denote some geodesic rays with $[c_n]=[q_n]$ and $[c]=[q]$ in $\partial \Gamma.$ We need to show that for each $t$ there exists an integer $k=k_t$ such that for all $n \geq k,$ we have $d(c_n(t), c(t)) \leq D,$ for some constant $D$ not depending on $c_n$ or $c.$

Using the stability of quasi-geodesics in hyperbolic spaces, there exists a constant $M=M(\cc_1,\cc_2,\delta) \geq 0$ such that the Hausdorff distance between $c_n$ and $q_n$ is bounded above by $M.$ Similarly, the Hausdorff distance between $c$ and $q$ is bounded above by $M.$ Let $t \in \mathbb{R}^{+}$, notice that using the above, there exists a point $q_t \in q$ such that $d(c(t),q_t) \leq M.$ Let $b_t$ be the furthest point in $b$ from $\go$ with $p(b_t)=q_t$. That is, $b_t $ is the unique point in $b$ with $p(b_t)=q_t$ and $d(b_t,\go) \geq d(x,\go)$ for any $x \in b$ with $p(x)=q_t.$ Let $s_t \in [0, \infty)$ with $b(s_t)=b_t$, since $b_n \rightarrow b,$ there exists some $k$ such that $d(b_n(s_t), b(s_t)) \leq 1$ for all $n \geq k.$ Since $p$ is coarsely distance decreasing, we have $d(p(b_n(s_t)), p(b(s_t)))=d(p(b_n(s_t)),p(b_t))=d(p(b_n(s_t)),q_t) \leq Kd(b_n(s_t),b_t))+C \leq K+C.$ Denote $p(b_n(s_t)) \in q_n$ by $x^t_n$. Notice that $d(x^t_n, q_t) \leq K+C$. By stability of quasi-geodesics, there exists a point $y^t_n \in c_n$ with $d(y^t_n,x^t_n) \leq M.$ Using the above, we get that $$d(c(t),y^t_n) \leq M+K+C+M=2M+K+C.$$ In particular, since $c_n,c$ are geodesics, if $s_n^t$ is the point in $[0, \infty)$ with $c_n(s_n^t)=y_n^t,$ then $s_n^t \in [t-(2M+K+C),t+(2M+K+C)]$. This implies that $d(c_n(t), c(t)) \leq D=2(2M+K+C)$ for all $n \geq k$ which finishes the proof. 
\end{proof}

\begin{corollary}
Fix a cocompact \CAT\text{ }cube complex $X$. If $L$ is the separation constant, then we have the following:

\begin{enumerate}
    \item There exists a coarsely distance decreasing map $v: (X,d) \rightarrow (\Gamma_L,d_{\Gamma_L})$.
    
    \item The map $v$ induces a continuous injection $i_v: \partial^{\kappa}_\infty X \rightarrow \partial \Gamma_L. $
    \item The map $v$ induces a continuous injection $i_v: \partial_{\kappa}X \rightarrow \partial \Gamma_L. $
\end{enumerate}

\begin{proof}
Part (1) of the statement is Corollary \ref{cor: new vertex map}. Part 2 follows using Lemma \ref{lem: injectivity of the map}, Lemma \ref{lem: continuous map gnrl settings},  Corollary \ref{cor: sublinearly Morse geodesics define boundary points} and part (1). Part(3) follows using Lemma \ref{lem:injection} and part (2).
\end{proof}

\end{corollary}

\section{ a remark on hierarchically hyperbolic groups}\label{sec: remark on hhg} 

The first part of Theorem \ref{introthm:cont injection} states that for a cocompact CAT(0) cube complex $X$ with a factor system, there exists a hyperbolic graph $\Gamma$ such that infinite $\kappa$-Morse geodesic rays projects to infinite unparamterized quasi-geodesics in $\Gamma$ defining a map $\partial_\kappa X \rightarrow \partial \Gamma$ where $\partial_\kappa X$ and $\partial \Gamma$ are the $\kappa$-Morse boundary and the Gromov boundary of $X$ and $\Gamma$ respectively.

It is natural to wonder if the same conclusion of Theorem \ref{introthm:cont injection} holds in the settings of hierarchically hyperbolic spaces.

\begin{question}
Let $\mathcal{X}$ be a geodesic hierarchically hyperbolic space and let $Y$ be a hyperbolic space such that Morse geodesic rays in $\mathcal{X}$ project to infinite unparameterized quasi-geodesics in $Y$. Do $\kappa$-Morse geodesic rays in $\mathcal{X}$ project to infinite unparameterized quasi-geodesics in $Y$?
\end{question}

We give an example showing that the answer to this question is no. In order to do that, we will use the following theorem by Rafi and Verberne (Theorem 1.1 in \cite{RafiVerberne}). For a surface $S,$ we let Map$(S)$ denote the mapping class group of $S$ and $\mathcal{C}(S)$ denote the curve graph of $S$.

\begin{theorem} \label{thm:badgeodesics} Let $S=S_{0,5}$ be the five-times punctured sphere and let $H=Map(S)$ be the mapping class group of $H$. There exists a finite generating set $A$ of $H$ and an infinite sequence of finite geodesics $b_n \in Cay(H,A)$ such that projections of $b_n$ to $\mathcal{C}(S)$ are not unparameterized $(K,C)$-quasi-geodesics for any $K,C>0.$
\end{theorem}

Now, let $H$ and $A$ be as in Theorem \ref{thm:badgeodesics}. Define $G=H \ast \mathbb{Z}$ and choose the finite generating set for $G$ by $A'=A \cup \langle a \rangle $, where $a$ is a generator for $\mathbb{Z}.$ Let $Y=\mathcal{C}(S) \ast \mathbb{R}$ be the hyperbolic space which is a ``tree of curve graphs". Notice that the action of $H$ on $\mathcal{C}(S)$ extends naturally to an action of $G$ on $Y$ giving a projection map from $G$ to $Y.$ Also, it's worth pointing out that the space $Y$ is precisely the maximal hyperbolic element of the HHS structure from \cite{Abbot}.

\begin{lemma}(Special case of Theorem 7.9 in \cite{Abbott2019} and Theorem 6.6 in \cite{Abbot})
If $\mathcal{X}=Cay(G,A')$, then $\mathcal{X}$ is a hierarchically hyperbolic space such that Morse geodesic ray in $\mathcal{X}$ project to infinite unparameterized quasi-geodesics in $Y.$
\end{lemma}

\begin{corollary}
Let $G,A',Y$ be as above and consider the hierarchically hyperbolic space $\mathcal{X}=Cay(G,A')$. There exists a $\sqrt{t}$-Morse geodesic ray in $\mathcal{X}$ whose projection to $Y$ is not an unparameterized $(K,C)$-quasi-geodesic for any $K,C>0.$
\end{corollary}

\begin{proof}

Recall that $G=H \ast \langle a \rangle $, were $\mathbb{Z}=\langle a \rangle.$ Notice that geodesics  $b_n$  from Theorem \ref{thm:badgeodesics} are all naturally geodesics in $\mathcal{X}=Cay(G,A')$, and hence the projection of the sequence $b_n$ to $Y$ are not unparameterized $(K,C)$-quasi-geodesics for any $K,C>0.$ We now use this sequence to build a $\kappa$-Morse geodesic $b$ ray with $\kappa=\sqrt{t}.$ More precisely, define $c_i$ by

$$c_1=\text{(first letter of $b_1$)}a\text{(first two letters of $b_1$)}a \dots \text{(all letters of $b_1$)}a.$$
$$c_2=\text{(first letter of $b_2$)}a\text{(first two letters of $b_2$)}a \dots \text{(all letters of $b_2$)}a.$$
$$\vdots$$
$$c_i=\text{(first letter of $b_i$)}a\text{(first two letters of $b_i$)}a \dots \text{(all letters of $b_i$)}a.$$
$$\vdots$$

We define $b$ as follows:

\begin{center}
    $b=c_1c_2c_3...c_i...$.
\end{center}

It is left as an exercise to the reader to check that $b$ is $\kappa$-Morse for $\kappa=\sqrt {t}.$

\end{proof}






\bibliographystyle{alpha}
\bibliography{bibliography}

\end{document}